\documentclass[10pt]{article}
\usepackage{latexsym,amsmath,amssymb,graphics}
\makeatletter
\@addtoreset{figure}{section}
\def\thefigure{\thesection.\@arabic\c@figure}
\def\fps@figure{h,t}
\@addtoreset{table}{bsection}

\def\thetable{\thesection.\@arabic\c@table}
\def\fps@table{h, t}
\@addtoreset{equation}{section}

\makeatother

\begin{document}

\newtheorem{theorem}{Theorem}[section]
\newtheorem{definition}[theorem]{Definition}
\newtheorem{lemma}[theorem]{Lemma}
\newtheorem{remark}[theorem]{Remark}
\newtheorem{proposition}[theorem]{Proposition}
\newtheorem{corollary}[theorem]{Corollary}
\newtheorem{example}[theorem]{Example}
\newtheorem{examples}[theorem]{Examples}

\newcommand{\bfi}{\bfseries\itshape}

\newsavebox{\savepar}
\newenvironment{boxit}{\begin{lrbox}{\savepar}
\begin{minipage}[b]{15.8cm}}{\end{minipage}
\end{lrbox}\fbox{\usebox{\savepar}}}

\makeatletter
\title{{\bf Asymptotic and Lyapunov stability of Poisson
equilibria}}
\author{Juan--Pablo Ortega$^1$,   V\'{\i}ctor Planas--Bielsa$^2$,
and Tudor S. Ratiu$^3$}
\addtocounter{footnote}{1}
\footnotetext{Centre National de la Recherche Scientifique,
D\'epartement de Math\'ematiques de Besan\c con, Universit\'e de
Franche-Comt\'e, UFR des Sciences et Techniques. 16, route de Gray.
F-25030 Besan\c con cedex. France. {\texttt
Juan-Pablo.Ortega@math.univ-fcomte.fr }}
\addtocounter{footnote}{1}
\footnotetext{Institut Non Lin\'eaire de Nice, UMR 129, CNRS-UNSA,
1361, route des Lucioles, 06560 Valbonne, France.
{\texttt Victor.Planas\_Bielsa@inln.cnrs.fr}}
\addtocounter{footnote}{1}
\footnotetext{Centre Bernoulli. \'Ecole Polytechnique F\'ed\'erale
de Lausanne. CH-1015 Lausanne.  Switzerland. {\texttt
Tudor.Ratiu@epfl.ch}}

\date{}
\makeatother
\maketitle

\addcontentsline{toc}{section}{Abstract}

\begin{abstract}
This paper includes results centered around three topics, all of
them related with the nonlinear stability of equilibria in Poisson
dynamical systems. Firstly, we prove an energy-Casimir type
sufficient condition for stability that uses functions that are not
necessarily conserved by the flow and that takes into account certain
asymptotically stable behavior that may occur in the Poisson
category. This method is adapted to Poisson systems obtained via 
a reduction procedure and we show in examples that the
kind of stability that we propose is appropriate when dealing with
the stability of the equilibria of some constrained systems.
Finally, we discuss two situations in which the use of continuous
Casimir functions in stability studies is equivalent to the
topological stability methods introduced by  Patrick {\it et al.}
~\cite{topological stability}. 
\end{abstract}

%\tableofcontents

\section{Introduction}

The use of the conserved quantities of a Hamiltonian flow in the
study of the stability of its solutions is a venerable topic that
goes back to Lagrange and Dirichlet. In the past decades these
ideas have been adapted to various setups: equilibria in
Poisson systems~\cite{eca, holm1, topological stability},
relative equilibria~\cite{patrick
92, singreleq, thesis, re, topological stability} and periodic and
relative periodic orbits~\cite{po, pos} of symmetric Hamiltonian 
systems, relative equilibria of symmetric Lagrangian
systems~\cite{slm}, and symmetric constrained systems~\cite{zenkov
energy momentum}, to list a few. All these results provide
sufficient conditions for the solution in question to be stable.

In this paper we will focus on the stability of the equilibria of 
Poisson dynamical systems. The main tools that one finds in the
literature concerning this case are the energy-Casimir method and
the topological stability methods introduced in~\cite{topological
stability}. The energy-Casimir method consists in finding a
combination of conserved quantities by the Hamiltonian flow,
typically the energy and the Casimir functions,
that  exhibits a critical point at the equilibrium with 
definite Hessian. Since the dynamics of the system is
confined to the level sets of this combination and, by the
Morse Lemma, in a coordinate chart about the equilibrium
these level sets are diffeomorphic to spheres
centered at the equilibrium,  stability follows. The
topological methods in~\cite{topological stability} rely on
a much more subtle confinement of the flow that takes
advantage not only of its conservation laws but also of
the  topological properties of the foliation of the
Poisson manifold by its symplectic leaves.

Energy
confinement is a very important tool in the symplectic
Hamiltonian context due to the absence of asymptotically stable
behavior. Energy methods are, to this day,   the only general
way to prove stability in more than two degrees of freedom. The
conservation of the phase space volume by the flow imposed by
Liouville's theorem does not necessarily hold in the
Poisson category. The first main result of this paper,
contained in Theorem~\ref{ultimate stability in Poisson},
adapts the standard energy-Casimir method to this
situation to allow the use of functions that are not
necessarily conserved by the flow but that can still be
used to conclude a certain kind of asymptotic
stability via the standard Lyapunov stability theorem.
This newly introduced notion of stability implies the
standard Lyapunov stability and will be referred to as
{\bfi weak asymptotic stability}. Theorem~\ref{ultimate
stability in Poisson} improves a previous version of the
energy-Casimir method (see~\cite{thesis} or  Corollary 4.11
in~\cite{pos}) where the conserved quantities 
confining the flow are also used to shrink the space on
which one checks the definiteness of the Hessian.
Theorem~\ref{ultimate stability in Poisson} shows that
\emph{any} conserved quantity can be used to shrink this space even
when that conserved quantity is not involved in the confinement of the
flow.

Theorem~\ref{subcasimir
stability} is the second main result of this paper. It 
adapts the stability condition in  Theorem~\ref{ultimate
stability in Poisson} to equilibria of Poisson systems obtained by
a certain reduction procedure that uses ideals in the Poisson
algebra of the functions on the manifold. Our interest is twofold.
First, there are some constrained mechanical systems that can be
described by reducing in this sense a bigger (unconstrained) system.
Moreover, the weakened kind of stability that Theorem~\ref{subcasimir
stability} allows us to conclude, coincides with the physically
relevant notion of stability in those situations, that is, the one
that describes the system when subjected to perturbations compatible
with the constraints. We illustrate  this point with a couple of
examples in Section~\ref{examples}: a light Chaplygin sleigh on a
cylinder and two coupled spinning wheels. Second, there are cases when
there are not enough conserved quantities to apply
Theorem~\ref{ultimate stability in Poisson} but, nevertheless, the
system can be reduced around the equilibrium and then the reduced
system has enough conserved quantities to use the theorem. 
Theorem~\ref{subcasimir stability} explains the meaning of having
this \emph{reduced} kind of stability. In particular, it shows the
role of sub-Casimir functions in stability computations.

The last section of the paper is dedicated to the study of the
relation between the topological stability methods
in~\cite{topological stability} with a generalized version of the
energy-Casimir method that we propose in the text based on the use
of local continuous Casimir functions of the Poisson manifold. To be
more explicit, the stability criteria in~\cite{topological
stability} are stated in terms of a set that, roughly speaking,
measures how far the space of symplectic leaves of a Poisson
manifold is from being a Hausdorff topological space. The general 
question that we  try to answer is under what circumstances this
set can be characterized as the intersection of level sets  of local
continuous Casimirs. Since this is not true in general, we provide two
sufficient conditions that are related to certain idempotency
of the set in~\cite{topological stability} and to the
possibility of separating regular symplectic leaves  by using
continuous Casimirs. The natural category where these
questions are posed is that of generalized foliated
manifolds; this is the context in which we have formulated
the main results in this section and where we have
obtained the Poisson case as a byproduct, considering it
as a manifold foliated by its symplectic leaves.

Before we start with the core of the paper we quickly review in a few
paragraphs the basic notions and terminology of Poisson
manifolds and generalized foliations that we will use throughout
the paper. The expert can safely skip the rest of this section.

\smallskip

\noindent {\bf Poisson systems.} Let $P$ be a smooth
manifold and let $C^{\infty}(P)$ be the algebra of smooth
functions on $P$. A Poisson structure on $P$ is a bilinear map $\{
\cdot , \cdot \} : C^{\infty}(P) \times C^{\infty}(P)
\longrightarrow C^{\infty}(P) $ that defines a Lie algebra
structure on $C^{\infty}(P)$ and that is a derivation on each
entry. The derivation property allows us to assign to each
function $F \in  C^{\infty}(P) $ a vector field $X _F \in
\mathfrak{X} (P)  $ via the equality
\[ X_H[F]:=\{F, H\} \text{ for every }F \in C^{\infty}(P).
\] The vector field $X_H \in \mathfrak{X}(P)$ is called
the {\bfi Hamiltonian vector field} associated to the {\bfi
Hamiltonian} function $H$.  
The derivation property of the Poisson bracket also implies that
for  any two functions $F,\,G\in C ^{\infty}(P)$, the value of the
bracket $\{F,\,G\}(z)$ at an arbitrary point $z\in P$ depends on
$F$ only through $\mathbf{d} F(z)$  which allows us to define  a
contravariant antisymmetric two-tensor $B\in\Lambda^2(P)$
by 
\[B(z)(\alpha_z,\,\beta_z)=\{F,\,G\}(z),\]
where $\mathbf{d} F(z)=\alpha_z \in T^\ast _z P$ and $\mathbf{d}
G(z)=\beta_z\in
T^\ast _z P$.  This tensor is called the {\bfi Poisson tensor\/} of
$M$.
The vector bundle map $B^\sharp:T^\ast  P\rightarrow TP$ naturally
associated to $B$ is defined by
$
B(z)(\alpha_z,\,\beta_z)=
\langle\alpha_z,\,B^\sharp(\beta_z)\rangle$.
Its range $\mathcal{E}:=B^\sharp(T^\ast  P) \subset TP$ is called
the {\bfi characteristic distribution} of the Poisson manifold
$(P, \{ \cdot , \cdot \})$. Its fiber at $z \in P $ is hence given by
$
\mathcal{E}_z=\lbrace X_H (z)\mid H \in
C^{\infty}(P)\rbrace
$. 
The distribution
$\mathcal{E} $ is a smooth generalized distribution which is always
integrable in the sense of Stefan~\cite{stefan,
stefan b} and Sussmann~\cite{sussmann}. Its maximal integral
submanifolds
$
\{\mathcal{L}\} $ are symplectic and are called the {\bfi  symplectic leaves} of
$(P, \{
\cdot , \cdot
\})$. The symplectic form $\omega _{\mathcal{L}}$ on the 
leaf $\mathcal{L}$ is uniquely characterized by the identity
\begin{equation*}
\omega _{ \mathcal{L}} (z)\left(  X_F(z) ,
X_G(z)\right):=\{F , G\}(z)
\quad \text{ for any } \quad F ,G \in C^{\infty}(P) \quad \text{ and for any }
\quad z \in \mathcal{L}.
\end{equation*} 
Since the symplectic leaves of $(P, \{
\cdot , \cdot
\})$ are the maximal integral leaves of a generalized
distribution, they form a {\bfi  generalized foliation} in the
sense of~\cite{dazord 1985}. This implies the
the existence of a chart  $(U, \varphi:U \rightarrow
\mathbb{R}  ^m)$ around any point
$z \in P $ such that if
$\mathcal{L}_z$ is the symplectic leaf containing $z$ then there
is a countable  subset $A \subset \Bbb R ^{m-n} $, with $m=\dim P
$ and
$n = \dim \mathcal{L}_z$, such that 
\begin{equation}
\label{expression dazord}
\varphi(U \cap \mathcal{L}_z)=\{ y \in \varphi
(U)\mid (y ^{n+1},
\ldots , y ^m)\in A\}.
\end{equation} 
Such a chart $(U,\varphi)$ is called a {\bfi foliation chart} for
the generalized symplectic foliation of $P$ around the point $z$. A
connected component of
$U \cap \mathcal{L}_z$ is called a {\bfi plaque} of the foliation
chart $(U, \varphi) $.  The point $z$ is said to be {\bfi  regular}
if the neighborhood $U$ can be shrunk so that all the leaves that
it intersects have all the same dimension. In that case, the plaques
coincide with the points of the form $(y ^1, \ldots, y ^n, y
^{n+1} _0, \ldots, y _0 ^m)\in \varphi (U) $ with $(y
^{n+1} _0, \ldots, y _0 ^m) $ constant. A leaf consisting of
regular points is said to be {\bfi  regular} and {\bfi  singular}
otherwise. The set of regular points of a generalized smooth
foliation is open and dense.

Some of the results proved in this paper will be first
given in the category of foliated manifolds. 
The corresponding results in the context of Poisson
manifolds are then obtained as corollaries.

\smallskip

\noindent {\bf Casimirs, local Casimirs, and first integrals of
foliations.} A function on a foliated manifold that is constant
on the leaves is called a {\bfi  first integral}
of the foliation. When we consider the particular case of a
Poisson manifold, the elements in
the center of the Poisson algebra
$(C^{\infty} (P),\{ \cdot , \cdot \})$, also called the {\bfi
Casimir} functions, are first integrals of the foliation of $P$
by its symplectic leaves. A {\bfi local Casimir} at
the point $z \in P$ is a function $C \in C^{\infty}(U_z)$ for some
open neighborhood $U_z \subset P$ of $z$ such that it is a Casimir of the
Poisson
manifold $(U_z, \{ \cdot , \cdot \}_{U_z})$ where the bracket $\{
\cdot , \cdot \}_{U_z} $ is the restriction of the bracket $\{
\cdot , \cdot \} $ on $P$ to $U_z$.

In general, non-trivial global Casimir functions may not exist. On the
other hand, local Casimirs are always available in the neighborhood
of a regular point. Indeed, if we think of the Poisson manifold $(P,
\{
\cdot , \cdot \})$ as a foliated space by its symplectic
leaves, the expression~(\ref{expression dazord}) 
allows us to find a chart $(U, \varphi:U \rightarrow
\mathbb{R}  ^m)$ around the regular point where
the  plaques of the symplectic foliation are the points of the form $(y ^1,
\ldots, y ^n, y
^{n+1} _0, \ldots, y _0 ^m)\in \varphi (U) $ with $(y
^{n+1} _0, \ldots, y _0 ^m) $ constant.
The functions that depend on
the last $m-n  $  coordinates are local Casimir functions of $(P, \{
\cdot , \cdot \})$ around $z$. 

\smallskip

\noindent {\bf Quasi-Poisson submanifolds and sub-Casimirs.} An 
embedded submanifold $S$ of $P$ which is Poisson in its
own right and is such that the inclusion
$i: S \hookrightarrow P $ is canonical is called a {\bfi
Poisson submanifold\/} of $P $. The Poisson structure on
$S $ is uniquely determined by the condition that the
inclusion be canonical, that is, there is no other Poisson
structure on $S $ relative to which the inclusion is
canonical. 

It turns out that in this paper we need a
slightly weaker condition. An embedded submanifold
$S$ of $P$ (without any Poisson structure on it) such that
$B^{\sharp}(s)\left(T_s^{\ast} P \right) \subset T_sS$ for
any $s \in S$ is called a {\bfi  quasi-Poisson
submanifold} of $P$. Every Poisson submanifold is
quasi-Poisson but the converse  is not true. 
As a corollary to  the main theorem
in~\cite{poissonreduction}, one can easily conclude that 
if $S $ is a quasi-Poisson submanifold of $P $, then there
is a unique Poisson structure $\{ \cdot ,\cdot\}_S$ on $S$
with respect to which the inclusion
$S\hookrightarrow P $  is a Poisson map, that is, there is
a unique induced Poisson structure on $S$ making it into
a Poisson submanifold of $P$. The Poisson bracket
$\{ \cdot , \cdot \}_S$ is  defined by $\{f,g\}_S(s) :=
\{F, G\}(s)$ where $F, G \in C^{\infty}(P)$ are arbitrary
local extensions of $f,g \in C^{\infty}(S)$ around the
point $s\in S$; this means that there is an open 
neighborhood $U  $ of $s$ in $P $ such that $f|_{S\cap
U}=F|_{S\cap U}$ and $g|_{S\cap U}=G|_{S\cap U}$. 

Thus, it is possible that the quasi-Poisson submanifold
$S$ of $P$ has its own Poisson structure (that is given
a priori) but it is not the one induced by the Poisson
structure of $P$. For a discussion of these issues see
\cite{ortega book}, sections 4.1.21 - 4.1.23.

Let $c \in C^{\infty}(S)$ be a Casimir function for the
Poisson manifold $(S, \{ \cdot , \cdot \}_S)$. Any  
extension $C \in C^{\infty}(P)$ of $c$ will be called a
{\bfi sub-Casimir} of $(C^{\infty} (P),\{ \cdot , \cdot
\})$.

Here is an example of the construction just described. Take
some Casimir functions $C_1, \ldots , C_k \in C^{\infty}(P)$ of 
$(P, \{ \cdot , \cdot \})$ and assume that a certain common 
level  set $S$ of these Casimirs is an embedded submanifold of
$P$. It is easy to check that 
$B^{\sharp}(s)\left(T_s^{\ast}P \right) \subset T_sS$ for
any $s \in S$ and hence $S$ carries a unique Poisson
bracket $\{\cdot , \cdot\}_S)$ such that 
$(S,\{\cdot ,\cdot\}_S)$ is a Poisson manifold with its
own Casimir functions that extend to sub-Casimirs on $P$.

\section{Stability in Poisson systems}

In this section we use some aspects of the geometry of  Poisson
manifolds to study the stability of the equilibria of 
Hamiltonian vector fields. 

Let $M$ be a manifold, $X \in \mathfrak{X}(M)$ a vector
field, $F_t$ the flow of $X$, and $m_e \in M $ an
equilibrium of $X $, that is, $X(m_e) = 0 $ or,
equivalently, $F_t(m_e) = m_e $ for all $t \in
\mathbb{R}$. Recall that $m_e$ is  {\bfi stable}, or
{\bfi Lyapunov stable}, if for any open neighborhood $U$
of $m_e$ in $M$ there is an open neighborhood
$V\subset U$ of $m_e$  such that $F_t(m) \in U$ for any
$m \in V$ and for any $t > 0$. The equilibrium $m_e$ 
is {\bfi  asymptotically stable} if there is a
neighborhood $V$  of $m_e$ such that  $F _t (V) \subset F
_s (V)$ whenever $t>s$  and $\lim\limits_{t \to \infty}F_t(V) =
m_e $, that is, for any neighborhood $W $ of $m_e $
there is a $T > 0 $ such that $F_t(V) \subset W $ if $t
\geq T $. If only the first condition holds, that is, $F _t (V)
\subset F _s (V)$ whenever $t>s$, we say that $m _e $ is {\bfi 
weakly asymptotically stable}. Note that 
\[
\text{{\bf asymptotic stability}}\Rightarrow \text{{\bf weak
asymptotic stability}}\Rightarrow\text{{\bf Lyapunov stability.}}
\]
Asymptotic stability cannot occur in
symplectic Hamiltonian systems due to Liouville's theorem ; only
Lyapunov stability is allowed. In the Poisson category, equilibria
lying in trivial symplectic leaves may be asymptotically stable.
However, if the symplectic leaf that contains the equilibrium is
at least two-dimensional, weak asymptotic stability is the
most we can hope for. 

The {\bfi linearization} of $X$ at the
equilibrium point $m_e$  is the linear map $L:T_{m_e} M
\rightarrow T_{m_e} M$  defined by $L (v) 
:=\left.\frac{d}{dt}\right|_{t=0} \left( T _{m_e} F _t
(v)\right)$ where $F_t$ is the flow of $X$  and
$v  \in T_{m_e} M$  is arbitrary. As is well known, the
study of the spectrum of the linear map $L$ gives
relevant information about the stability of the
equilibrium $m_e$.  The equilibrium $m_e \in
M$ is {\bfi linearly stable} (respectively
{\bfi unstable)} if  the origin is a stable (respectively
unstable) equilibrium for the linear dynamical system on
$T_{m_e} M$  defined by $L$. The equilibrium
$m_e$ is {\bfi spectrally stable}    (respectively
{\bfi unstable)} if the spectrum  of the linear map $L$
lies in the (strict) left-half plane or on the imaginary
axis (respectively at least one eigenvalue has strictly
positive real part). Lyapunov and linear stability imply
spectral stability. If all the eigenvalues of $L$ have
strictly negative real part, that is, they lie in the
(strict) left-half plane, the system is
asymptotically stable.

\subsection{Linearization of Poisson dynamical systems and
linear stability}

Consider
a Hamiltonian vector field $X_H$ on the Poisson manifold $(P,\{
\cdot , \cdot \})$, let $z_e \in P$ be an equilibrium of
$X_H $, and $L: T_{z_e} P \rightarrow T_{z_e} P$
the linearization of $X_H $ at $z_e $. If $z_e $ is
regular (in particular, when  $P$ is a symplectic
manifold) there are restrictions on the eigenvalues of
$L$ that do not allow us to  conclude the Lyapunov
stability of $z_e$ from its spectral 
stability (see, for instance, Theorem 3.1.17 in~\cite{fom}).  As
will be shown below, this restriction disappears, in general, for 
equilibria lying on singular symplectic leaves.

In order to present the following lemma, whose proof is a
straightforward computation, we recall that there exists a chart 
$(U,
\varphi)$ around any point $z \in P $ in the $2n+r$ dimensional
Poisson manifold $(P, \{
\cdot ,
\cdot
\})$   such that $\varphi (z) = {\bf 0} $ and that the associated
local coordinates, denoted by $(q ^1,\ldots, q ^n, p _1,
\ldots, p _n,z ^1, \ldots, z ^{r})$, satisfy
$
\{q ^i, q ^j\}=\{ p _i, p _j\}=\{q ^i, z ^k\}=\{p _i, z
^k\}=0$  and  $\{q^i, p _j\}= \delta^i_j$,
for all $i,j,k $ such that  $1 \leq i,j\leq n $, $1\leq k\leq r$. 
For all such that $k,l $, $1\leq k,l\leq r $, the Poisson bracket 
$\{z ^k, z^l\}$ is a function of the local coordinates $z
^1, \ldots, z ^{r} $ exclusively and vanishes at $z$. Hence,
the restriction of the bracket $\{
\cdot , \cdot \} $ to the coordinates $z^1, \ldots,
z^{r}$ induces a Poisson structure on an open
neighborhood $V$ of the origin in $\mathbb{R}^r$ whose
Poisson tensor will be denoted by $\mathcal{R} \in
\Lambda ^2 (V)$. This Poisson structure on $V $ is
called  the {\bfi  transverse Poisson structure}  of $(P,
\{ \cdot , \cdot\})$ at $z$ and is unique up to Poisson
isomorphisms. The coordinates of the local chart that we
just described are called  {\bfi  Darboux-Weinstein}
coordinates~\cite{weinstein83}.

\begin{lemma}
\label{linearization darboux weinstein} 
Let $z_e$ be an equilibrium of the Hamiltonian dynamical
system on the Poisson manifold $(P,\{ \cdot , \cdot \})$
and let
$(\mathbf{q}, {\bf p}, {\bf z}) $ be a Darboux-Weinstein chart
around $z$. Denote by $\mathbf{x}:=(\mathbf{q}, {\bf p})$ and by
$J $ the $n \times  n $ square matrix given by
\[
J=\left(
\begin{array}{cc}
0 &I _n\\
-I _n& 0
\end{array}\right).
\]
The linearization $L$ of $X _H $
at the equilibrium $z_e$ in the coordinates
$({\bf x}, {\bf z})$  takes the form
\begin{equation}
\label{linearization of L}
L=\begin{pmatrix} \mathcal{S} & \mathcal{Q}\\ 0 &
\mathcal{P}
\end{pmatrix}
\end{equation}
where 
\[
\mathcal{S} ^i _j= \sum _{p=1}^{2n}J ^{i p} \frac{\partial
^2H}{\partial x ^p\partial x ^j}({\bf 0}, {\bf 0}), 
\quad
\mathcal{P} ^k _l=\sum _{p=1} ^r \frac{\partial \mathcal{R}^{k
p}}{\partial z ^l}({\bf 0})\frac{\partial H}{\partial z ^p}({\bf
0}, {\bf 0}), 
\quad \text{and}\quad
{\cal Q}^i _l= \sum _{p=1}^{2n} J ^{ i
p}\frac{\partial ^2 H}{\partial x ^p\partial z^l}
({\bf 0}, {\bf 0}).
\]
\end{lemma}

\noindent\textbf{Proof.\ \ } The result is obtained by
differentiating the expression of the Hamiltonian vector field at the
equilibrium in Darboux-Weinstein coordinates and by taking into
account that the matrix 
$J$ is constant, that
$\mathcal{R}({\bf 0})$ is zero, and that $\mathcal{R}$
depends only on the
${\bf z}$ variables. \quad
$\blacksquare$

\medskip

We now use~(\ref{linearization of L}) to give a
characterization of the structure of the eigenvalues of the
linearized vector field $L$ in the Poisson context. The proof of the
following proposition is a straightforward computation.

\begin{proposition} 
\label{eigenvalues and so on}
In the situation described in the previous
lemma denote by 
$\{
\lambda_1 ,
\ldots ,
\lambda_{2n}\}$ the eigenvalues of the infinitesimally
symplectic matrix $\mathcal{S}$, counted with their
multiplicities, and let 
$\{u_1 ,\ldots , u_{2n}\}$ be a basis of corresponding
eigenvectors. Assume that the matrix
$\mathcal{P} $ is diagonalizable, let
$\{\mu_1,\ldots ,\mu_r\} $ be its eigenvalues counted
with their multiplicities, and 
$\{v_1 , \ldots , v_r \}$ a basis of eigenvectors. Then
the matrix $L$ has eigenvalues $\{ \lambda_1 , \ldots,
\lambda_{2n}, \mu_1 , \ldots , \mu_r\}$.  If for any
eigenvalue $\mu_j$ we have that $(\mathcal{S}-\mu_j
I)^{-1}\mathcal{Q}v_j$ is not empty then $L$ is
diagonalizable with corresponding  basis of eigenvectors 
\[\{(u_1,0),\ldots (u_{2n},0),
(-w_1,v_1), \ldots,  (-w_r,v_r)
\},
\]
where $w_j \in (\mathcal{S}-\mu_j I)^{-1}\mathcal{Q}v_j$,
$j = 1, \dots, r$ are arbitrary but subjected to the
condition that if $v_j = v_k $ then $(w_j,v _j)$ and $(w _k, v _k)$
are chosen to be linearly independent.
\end{proposition}

The eigenvalues $\{ \lambda_1, \ldots, \lambda _{2n}\}$ 
satisfy the symplectic eigenvalue theorem since
$\mathcal{S}$ is infinitesimally symplectic. However,
the eigenvalues $\{\mu_1, \ldots, \mu _r\}$ may lie, in
principle, anywhere in the complex plane. Hence  Poisson
dynamical systems may exhibit asymptotic behavior. There
are three specific situations that should be singled out:
\begin{itemize}
\item  None of the eigenvalues of $\mathcal{P}$ coincides with one
of the  eigenvalue of $\mathcal{S}$. In this case the matrices 
$(\mathcal{S}-\mu_j I)$, $1 \leq j\leq r$, are invertible
and the whole linear system $L$ is diagonalizable. 
\item $\mu_i=\lambda_j$ for some $i , j$ but
$(\mathcal{S}-\mu_i I)^{-1}\mathcal{Q}v_i$ is not empty. Then there is a
passing of eigenvalues but they do not interact in the sense
that they correspond to different blocks in the linearized
system. We will call this
situation  {\bfi uncoupled 
passing}. 
\item If in the previous case $(\mathcal{S}-\mu_i I)^{-1}\mathcal{Q}v_i$ is
empty then the linear system is not  diagonalizable anymore and the
passing of eigenvalues mixes blocks of the infinitesimally
symplectic part and the transversal one. We will call this
situation  {\bfi coupled 
passing}.
\end{itemize}
With
these remarks in mind, we get the following.

\begin{proposition} 
Let $(P, \{ \cdot , \cdot \},H)$ be a Poisson
dynamical system and $z_e \in P$ an equilibrium point of
$X _H $. If the linearization $L$ of $X _H  $ at $z_e$
exhibits a coupled passing then the system is linearly
unstable.
\end{proposition}
\textbf{Proof}
\noindent The existence of a coupled passing implies the occurrence
in $L$ of a nondiagonal block in its Jordan canonical form. The
flow of the linear dynamical system induced by $L$, when restricted
to the space generated by the associated Jordan basis, exhibits an
unstable behavior and the result follows.
$\blacksquare$

\begin{corollary} 
\label{regular linearization}
Consider the linearization $L$ of a Poisson
dynamical system $(P,\{ \cdot , \cdot \}, H)$ around an
equilibrium $z_e  \in  P  $  lying on a regular
symplectic leaf $\mathcal{L}$.  Let $\{ \lambda_1,
\ldots , \lambda_{2n}\}$ be the eigenvalues of the
infinitesimally symplectic block $\mathcal{S}$. Then
\begin{description}
\item [(i)]  
$\mathcal{P}=0$.
\item [(ii)]  The vectors $u \in T_{z_e}P$ that satisfy $Lu=
\lambda u $ for some $\lambda\neq 0 $  lie in
$T_{z_e}\mathcal{L}$. In particular, the unstable directions of $L$
are tangent to the symplectic leaf of  $P$ that contains the
equilibrium.
\item [(iii)]  If $S^{-1} \mathcal{Q}v_j$ is not empty for any
$v _j$ as in Proposition~\ref{eigenvalues and so on}   then 
$0$ is the only eigenvalue in addition to  $\{
\lambda_1, \ldots ,
\lambda_{2n}\}$.
\end{description}
\end{corollary}

\noindent\textbf{Proof.\ \ }The first part follows from the expression
for $\mathcal{P}$ provided in Lemma~\ref{linearization darboux
weinstein} and from the fact that
$\mathcal{R}=0$ in an open neighborhood of
$z_e $ that contains only regular points. The unstable
directions are the vectors in the eigenspaces
corresponding to strictly positive eigenvalues. Then the
points {\bf (ii)} and {\bf (iii)} follow from the
expression of
$L$ in Lemma~\ref{linearization darboux
weinstein} using that on the set of regular points
$\mathcal{R}=\mathcal{P}=0$.
\quad $\blacksquare$

\subsection{Non linear stability in Poisson dynamical systems}
\label{stability Poisson section} 
As noted in the previous subsection, the array of linear
tools available to conclude nonlinear stability of 
equilibria  of a Poisson dynamical system is very limited.  In this
section we will state a result that provides a sufficient condition
for such  equilibria to be Lyapunov or weakly asymptotically stable.
This result is based on the use of local Casimirs and conserved
quantities of the dynamical system in question and is related to the
classical energetics methods (also called Dirichlet criteria)
in~\cite{eca, topological stability}. Our approach builds on an
improvement of the classical result in~\cite{eca} that was carried out
in~\cite{thesis}  (see Corollary 4.11 in~\cite{pos}).

The proof of our main result will be based on a classical result 
of Lyapunov that states that \textit{if $m_e \in  M $ is
an equilibrium of the vector field $X \in
\mathfrak{X}(M) 
$ with flow $F _t  $  and  there exists a positive
function $L\in C^{\infty}(U)$ around $m_e$, with $U$ an
open neighborhood of $m_e$, such that $\dot{L}(m):=
\left.\frac{d}{dt}\right|_{t=0}L(F _t (m))\leq 0$, for
any $m \in U \setminus \{m_e\}$, then $m_e$ is a Lyapunov
stable equilibrium}. We recall that a function $f \in
C^{\infty}(M)$ is said to be {\bfi  positive} around
$m_e\in M$ if $f (m_e)=0$ and there is an open
neighborhood $U_{m_e}$ of $m_e$ such that $f(m)> 0$, for
all $m \in  U _{m_e}\setminus\{m_e\} $. \textit{If
$\dot{L}(m) <0$ for all $m \in 
U_{m_e}\setminus\{m_e\}$, then $m_e$ is asymptotically
stable}. See e.g. Theorem 1, Chapter 9, \S3 in
\cite{HiSm1974} for a proof of these statements; the
infinite dimensional versions of these assertions can be
found in Theorems 4.3.11 and 4.3.12 of \cite{mta}. Any
positive function $L$ in the statement of Lyapunov's
theorem is usually called a {\bfi Lyapunov
function}. Its construction for specific dynamical
systems is by itself a very active research subject. 

In the case of Hamiltonian mechanics, the
Hamiltonian and the Casimirs of the Poisson phase space are natural candidates
to be used in Lyapunov's theorem.
If, additionally, the system has a symmetry to which one
can associate a momentum  map, its components are
conserved quantities that sometimes can be used for the
same purpose. The use of all conserved quantities of a
dynamical system in the study of the stability of
equilibria to form Lyapunov functions is known under the
name of energy-momentum methods.  However, it should be
noted that, apart from conserved quantities, Lyapunov's
theorem can be applied with the more general class of
functions whose time derivative is strictly negative.
The existence of these functions implies the asymptotic
stability of the equilibrium in question. In the
symplectic context  this is impossible. This behavior,
allowed only for Poisson Hamiltonian systems, is used 
in the main theorem of this subsection.

In the sequel we will use the following notation. Let
$P$ be a smooth manifold, $f  \in C^{\infty} (P) $ a smooth
function,  $z_e \in P$ a critical point of $f$  (that 
is, $\mathbf{d} f (z_e)=0 $), and $U$ an open
neighborhood of $z_e$. The {\bfi Hessian} of $f$ at the
critical point $z_e$ is the symmetric bilinear form
$\mathbf{d}^2f (z_e):T_{z_e} P \times T_{z_e}
P\rightarrow \mathbb{R}$ given by
$\mathbf{d}^2 f(z_e)(v,w):=v\lbrack W \lbrack f
\rbrack\rbrack$, where $v,w \in T_{z_e} P$ and $W \in
\mathfrak{X}(U)$ is an arbitrary extension of $w$ to a
vector field on $U$. The fact that $z_e$ is a critical
point of $f$ ensures that this definition is independent
of the extension
$W$ of $w$.

\begin{theorem}
\label{ultimate stability in Poisson}
Let $(P , \{ \cdot , \cdot \} ,H)$ be a Poisson dynamical 
system, $z_e$ an equilibrium point of the
Hamiltonian vector field $X_H$, and $C_1, \ldots , C_k : P \to
\mathbb{R}$ conserved quantities,
that is, $\{C _i, H \} =0 $, $i \in \{1, \ldots,
k\}$. Let $F : P \rightarrow  \mathbb{R} $ be a function
satisfying at least one of the following two conditions
\begin{description}
\item [(i)] $\{(F-F (z_e))^2,H  \}(z)\leq 0$ for all $z \in P
\setminus
\{z_e\}$ and $\{F-F (z_e),H  \}(y)\leq 0$  for all the points $y\in P
\setminus
\{z_e\}$ such that 
$\{(F-F (z_e))^2,H  \}(y)=0$. 
\item  [(ii)] $\{(F-F (z_e))^2,H  \}(z)\leq 0$ for all $z \in P
\setminus
\{z_e\}$ and $F-F (z_e) $ is positive on $P $  around $z _e
$.
\end{description}
Assume that there exist constants $\{\lambda_0,
\lambda_1, \ldots, \lambda_k, \mu\}$ such that
$$
\mathbf{d}(\lambda_0 H   + \lambda_1 C_1 + \ldots + 
\lambda_k C_k+\mu F)(z_e)=0
$$
and the quadratic form 
\begin{equation}
\label{definiteness hypothesis}
\mathbf{d}^2(\lambda_0 H + \lambda_1 C_1 + \ldots
+ \lambda_k C_k+\mu F)|_{W \times W}(z_e) 
\end{equation}
is positive definite, where 
\[
W=\ker \mathbf{d}H   (z_e)\cap \ker
\mathbf{d}C_1(z_e)\cap \ldots\cap\ker \mathbf{d}C
_k(z_e).
\]  
Then $z_e$ is a weakly asymptotically stable equilibrium (and
hence Lyapunov stable). If the inequality  $\{(F-F (z_e))^2,H
 \}(z)\leq 0$ is strict for every $z \in P\setminus\{z_e\}$ then
$z_e$ is asymptotically stable (this can only happen if the
symplectic leaf that contains the equilibrium is trivial).
\end{theorem}
\textbf{Proof}
Consider the functions $l _1, l _2\in C^{\infty}(P)$ 
defined by 
\begin{align*}
l_1(z)&:=\sum_{j=0}^{k}\left( \lambda_j C_j(z) + \mu
F(z)\right)-\left(\lambda_j C_j(z_e) + \mu F(z_e)\right),
\\
l _2 (z)&:=\sum_{j=0}^k
\frac{1}{2}\left((C_j(z)-C_j(z_e))^2 + (F(z) -
F(z_e))^2\right),
\end{align*}
where we denote $C _0:=H $.
Notice that $l_1(z_e)=0$ and that, by hypothesis,
$\mathbf{d}l_1(z_e)=0$ which implies that
$\mathbf{d}^2 l_1(z_e)$ is well defined.  Moreover, 
hypothesis~(\ref{definiteness hypothesis}) is
equivalent to the statement that 
$\mathbf{d}^2 l_1(z_e)|_{W \times W}$ is positive
definite. Additionally, $l_2(z_e) = 0$, $\mathbf{d}
l_2(z_e)=0$, and hence $\mathbf{d}^2 l_2(z_e)$ is  well
defined. A straightforward computation shows that
$\mathbf{d}^2 l_2(z_e)$ is positive semidefinite
with  kernel equal to the space $W$. A result due to
Patrick (see~\cite{patrick 92}) shows that in these
circumstances there exists a constant $r > 0$ such that
for any $\epsilon \in (0,r]$ the Hessian $\mathbf{d}^{2}
(\epsilon l _1+ l _2)(z_e)$ is positive definite. 

Let $L _\epsilon:=\epsilon l _1+ l _2 $. The positive
definiteness of $\mathbf{d} ^2L _\epsilon (z_e)$
implies  that $L _\epsilon$ is a positive function on
an open neighborhood $U$ of $z_e$ whose level sets are, by the Morse
lemma, diffeomorphic to concentric spheres centered at the
equilibrium $z _e$. Additionally, both hypotheses {\bf (i)} or
{\bf (ii)} imply that the constant $\epsilon $ can be chosen small
enough so that the time derivative
\begin{equation}
\label{computation of dot}
\dot{L}_\epsilon(z)=\{L
_\epsilon, H  \}(z)
=\frac{1}{2}\{(F-F (z_e))^2,H\}(z)+ \epsilon \mu\{
F-F (z_e),H\}(z)\leq 0,
\end{equation}
for any $z \in  P
\setminus \{z_e\}$. This implies that if $F _t $ is the flow of $X _H
$, the basis of open neighborhoods of $z _e $ given by the sets  $U
_\lambda:= L _\epsilon ^{-1}([0,
\lambda))$, with $\lambda $ small enough, satisfies $F _t(U _\lambda)
\subseteq F _s(U _\lambda) $, provided that $t \geq s $. This proves
the weak asymptotic stability of $z  _e$. 

If $\{(F-F (z_e))^2,H\}(z)< 0$  for every $z \in P\setminus\{z_e\}$
then $\epsilon $ can be chosen so that the positive function $L
_\epsilon $ is such that
$\dot{L}_\epsilon(z)<0 $ for any $z \in  P
\setminus \{z_e\}$ (see~(\ref{computation of dot})). Lyapunov's
theorem proves the asymptotic stability of $z _e$.
$\blacksquare$

\begin{remark}
\label{differences with previous}
\normalfont
The main differences between this result and those
already  existing in the literature are:
\begin{description}
\item [(i)] It takes advantage of the possible existence
of strict Lyapunov functions and hence is
capable of obtaining the Lyapunov stability of an
equilibrium as a corollary of an asymptotically stable
behavior. This feature allows us to prove stability in
some examples where no other available energy method is
applicable. 

In order to illustrate this
point consider the following example. The two dimensional
Toda lattice admits a  Poisson formulation~\cite{Bloch
Toda} by taking the bracket
$\{x,y\}=-x$ and the Hamiltonian function
$H(x,y)=x^2 + y^2$. This system has an equilibrium point at
$z_e=(0,b)$ for any $b \in \mathbb{R}$. The equilibrium
$(0,0)$  is obviously Lyapunov stable since
$\mathbf{d}H(0,0)=0$ and $\mathbf{d}^2H(0,0)>0$. The  
equilibria of the form $z_e=(0,b)$ with $b>0$ are
weakly asymptotically stable. This can be proved using the
previous theorem by taking the Hamiltonian as conserved
quantity and the function $F(x,y):=x$. The function $F $
satisfies the hypothesis {\bf (i)} in Theorem~\ref{ultimate stability
in Poisson} since
$\{F^2,H\}=-4 x^2 y\leq 0$  and $\{F,H\}=-2 x y=0 $  when
$\{F^2,H\}=0 $, in an open neighborhood of 
$z_e=(0,b)$ with $b>0$. If
$b<0$ the equilibrium is unstable since the linearization has an
eigenvalue with positive real part.
We emphasize that the stability of the points in the case $b>0 $ 
are uniquely due to their weak asymptotically stable
behavior. 

\item [(ii)] Unlike the approach taken in the treatment
of many standard examples (see for instance~\cite{ims})
this theorem shows that one  does not need to take
arbitrary functions of the conserved quantities in the
expression~(\ref{definiteness hypothesis}).  Indeed,
only linear combinations are needed. This is a
consequence of the fact that the form whose definiteness
needs to be studied is restricted to the space $W$.
\item [(iii)] Since the constants $\{\lambda_0,
\lambda_1, \ldots, \lambda_k, \mu\}$ are allowed to be
zero we have the freedom not to use a local conserved
quantity in the definiteness condition~(\ref{definiteness
hypothesis}) but to still take advantage of its
existence to shrink the space $W$. This is an improvement with
respect to the results in~\cite{thesis}  (see Corollary 4.11
in~\cite{pos})

In order to visualize this better
consider the following example. Let $(\mathbb{R} ^3,\{\cdot
, \cdot
\},H)$ be the Poisson dynamical system whose Poisson bracket is
given by the Poisson tensor that in Euclidean coordinates takes
the form
\[ B(x,y,z)=\begin{pmatrix} 0 & 0 & y\\ 0 & 0 & -x\\ -y & x
& 0
\end{pmatrix}
\]
and where  $H(x,y,z)=a z$ with $a \in \mathbb{R}$ a
nonzero constant. The function
$C(x,y,z)=\frac{1}{2}\left(x^2 + y^2 \right)$ is a
Casimir for this Poisson structure and  every point of
the form $(0,0,z_0)$ is an equilibrium of  $X_H$. Note
that $\mathbf{d}(H-\lambda C)(0,0,z_0) \not= 0$
for any $\lambda \in  \mathbb{R}$. Nevertheless, we can
still apply the previous theorem to conclude the
Lyapunov stability of $(0,0,z_0)$ by taking the
combination $\lambda_0 H + \lambda_1 C $  with
$\lambda _0=0 $ and $\lambda_1=1   $. With these choices,
$W={\rm ker}\left(\mathbf{d}H(0,0,z_0)\right)$ and
$\mathbf{d}^2 C(0,0,z_0)|_{W \times W}$ is positive
definite. The stability of these equilibria can also be handled using
the topological methods in~\cite{topological stability}. 
\end{description}
\end{remark}

\begin{remark}
\normalfont
The most efficient way to apply Theorem~\ref{ultimate stability in
Poisson} in order to establish the stability of  a given equilibrium
consists of looking at the Hamiltonian system obtained by restriction
of the original one to an arbitrarily small neighborhood of the
equilibrium. The advantages of proceeding in this way are based on
the fact that the restricted system has, in general, more conserved
quantities than the original one. We illustrate this remark with the
following specific example.

Consider the manifold
$P:=\mathbb{T}^2\times  \mathbb{R} $ endowed with the Poisson
structure given by the tensor that in coordinates $(\theta,
\varphi,x) $ is expressed as
\[
B(\theta,\varphi,x)=
\begin{pmatrix}
0 & 0 &1\\
0 &0 & - \alpha\\
-1 &\alpha &0
\end{pmatrix},
\qquad \alpha \in  \mathbb{R}\setminus \mathbb{Q}.
\]
Let  $H \in  C^{\infty}(P)  $  be the function defined by
$H(\theta,\varphi,x):=x ^2-\cos \theta $. The associated
Hamiltonian vector field $X_H$ has an equilibrium at
the point $z_e:=(0,0,0) $ whose stability we show using 
Theorem~\ref{ultimate stability in Poisson}. Even though the
Poisson manifold $P$ has no globally defined Casimir functions,
any locally defined function of the form $C= \alpha \theta +
\varphi $ is a local Casimir. We can use this local Casimir to
establish the Lyapunov stability of $z_e$. Indeed,
$\mathbf{d}H (z_e)=0 $ and $\mathbf{d} ^2H (z_e)|
_{W\times  W}>0 $, with $W=\ker
\mathbf{d} C (z_e) $. In section~\ref{Two coupled
spinning wheels} we will describe a mechanical system that is
closely related to this example.

\end{remark}

\begin{remark}
\label{lyapunov function and negative eigenvalues}
\normalfont
In most applications, the
conserved quantities in the statement of the theorem are
local Casimir functions, components of momentum maps, and the
Hamiltonian. 
A good way to find the functions 
$F$ is to look for purely negative eigenvalues of
the linearization of $X_H$  at the
equilibrium $z_e $ that do not have a positive
counterpart, as will be shown below. Notice that by
Corollary~\ref{regular linearization} this is only
possible when the equilibrium $z_e$ is lying on a
singular symplectic leaf of the Poisson manifold. More
explicitly, suppose that the linearization has such a
negative eigenvalue
$-\lambda$ with eigenvector $v$. Take local coordinates  
$(y_1,\ldots , y_n)$ such
that $v=\frac{\partial}{\partial y_n}$. Since  the
function $F_v(y_1,\ldots , y_n):=y_n$ satisfies
$\{y_n^2,H\}=X_H[y_n^2]=2 y_n \dot{y}_n =-2 \lambda y_n^2
+ {\rm h.o.t.}$, it is a good candidate to be used as the
function $F$ in the statement of the theorem. This
procedure has been used in the first example in
Remark~\ref{differences with previous}. 
\end{remark}

\subsection{Ideal reduction and ideal stability for Poisson systems}

We start this section by describing new Poisson structures on
some submanifolds of a Poisson manifold that can be obtained by
looking at the ideals  of its Poisson algebra of smooth
functions. We will refer to the construction that
will be presented as {\bfi   ideal reduction} for it is a
particular case of the Poisson reduction procedures
in~\cite{poissonreduction, poisson reduction singular, ortega
book}. 

This reduction technique is used later in this section to define
a weaker notion of stability, called  ${\cal I}${\bfi
-stability}, and to establish a
sufficient condition for it to hold. As the examples in
the next section show, the use of
${\cal I}$-stability is a very sensible way to deal with
the physically relevant stability properties of
equilibria in Hamiltonian systems subjected to
semi-holonomic constraints.

Let $P$ be a smooth manifold and $\mathcal{F}\subset 
C^{\infty}(P)$ be a family of smooth functions. Denote by
$\mathcal{V} _{\mathcal{F}} \subset P$ the {\bfi 
vanishing subset} of $\mathcal{F}$, defined as the 
intersection of the zero level sets of all the elements
of ${\cal F}$. For a subset $S \subset P$ define
its {\bfi  vanishing ideal} ${\cal I}(S)$ as the set of 
functions $f \in C^{\infty}(P) $ such that $f(S)=\{0\}$.
Notice that ${\cal I}(S)$ is obviously an ideal of
$C^{\infty}(P) $ with respect to the standard multiplication of
functions.  Notice also that for every subset $S \subset
P$ and for every ideal $\mathcal{J} \in C^{\infty}(P)$
we have $S \subset \mathcal{V}_{\mathcal{I}(S)} $ and
$\mathcal{J} \subset \mathcal{I}\left(\mathcal{V_J}
\right)$. These inclusions are in general strict.
However, if $S$ is a closed embedded submanifold of
$P$ then the first inclusion is actually an equality due to the
smooth version of Urysohn's lemma. Moreover, in this
particular case, the quotient algebra
$C^{\infty}(P)/\mathcal{I}(S)$ can be identified with
$C^{\infty}(S)$, the algebra of smooth functions on $S$
with respect to its own smooth manifold structure, via
the map that assigns to any $f  \in C^{\infty}(S) $ the
element $\pi(F) \in C^{\infty}(P)/\mathcal{I}(S)$, where
$F \in  C^{\infty} (P) $ is an arbitrary extension of
$f$ and $\pi:C^{\infty}(P) \rightarrow
C^{\infty}(P)/\mathcal{I}(S)$ is the projection.  We will
say that an ideal $\mathcal{I}\subset C^{\infty}(P)$ is 
{\bfi regular} if its vanishing set
$\mathcal{V_I} \subset P$ is a closed and embedded submanifold of
$P$. 

In the sequel we will focus our attention on finitely 
generated Poisson ideals. Let $(P, \{ \cdot , \cdot \})$
be a Poisson manifold and $\mathcal{F}=\{f_1, \ldots,
f_n\} \subset  C^{\infty}(P)$ be a finite family of
elements in $C^{\infty} (P) $. We will say that
$\mathcal{F} $ generates a {\bfi   Poisson ideal}
 if for any function $f 
\in  C^{\infty} (P)$ and any $i \in \{1, \ldots, n\} $
there exist functions $\{ h_{i1}, \ldots, h_{in}\}
\subset  C^{\infty}(P)$ such that
\[
\{f, f_i\}=\sum _{j=1}^{n}h _{ij} f _j.
\]
Denoting 
\[
\mathcal{I} ( \mathcal{F}) : = \left\{\sum_{k=1}^n
g_k f_k \,\Big|\, g_k \in C ^{\infty}(P) \right\}
\]
note that the  condition above is equivalent to the
statement that $\mathcal{I}(\mathcal{F}) $ is an ideal
in the Poisson algebra $C ^{\infty}(P) $, that is, it is
an ideal relative to both the usual multiplication of
functions as well as the the Lie bracket $ \{\cdot,
\cdot\}$. Note that \textit{if the vanishing
subset $\mathcal{V}_{\mathcal{F}}$ of  $\mathcal{F} $ is
an embedded submanifold of $P$ then $\mathcal{V}
_{\mathcal{F}}$ is a quasi-Poisson submanifold of $P$}.
Indeed, for any $f \in C^{\infty} (P) $, $f _i \in
\mathcal{F} $, and $z  \in \mathcal{V} _{\mathcal{F}}$,
there exist functions $\{h _1, \ldots, h _n\}
\subset  C^{\infty}(P)$ such that
\[
\langle \mathbf{d} f  _i (z),  X _f (z)\rangle 
= \{f_i, f\}(z)=\sum _{i=1}^{n}h _i(z) f _i(z)
=0,
\]
which shows that $B ^\sharp (z) (T ^\ast_z P) \subset 
T _z \mathcal{V}_{F}$, as required. Since the embedded
submanifold $\mathcal{V} _{\mathcal{F}} $ is
quasi-Poisson, it has a Poisson bracket $\{ \cdot , \cdot
\}_{\mathcal{V} _F}$ given by
$\{f, g\}_{\mathcal{V} _F}(z):=\{F , G\}(z)$, where $F, G
\in C^{\infty}(P)$ are arbitrary local extensions of $f,g 
\in C^{\infty}(\mathcal{V} _F)$ around the point $z\in
\mathcal{V} _F$. We recall that the extensions to $P$ of
the Casimir functions of $(\mathcal{V} _F,\{ \cdot ,
\cdot \}_{\mathcal{V} _F})$ are called
sub-Casimirs of $(P,\{ \cdot , \cdot \})$.

The construction that we just carried out can be locally reversed,
that is, given an injectively immersed quasi-Poisson submanifold $S$
of  $(P,\{ \cdot , \cdot \})$ any point $z \in S $ has an open
neighborhood $V _z  $ of $z  $ in $S$ such that the vanishing ideal
${\cal I}(V _z )$ is a Poisson ideal generated by a finite family
of smooth functions on $P$ with ${\rm codim}\, S $ elements.
Indeed, choose $V _z $ small enough so that it is an embedded
submanifold of $P $ and that, at the same time, is
contained in  the domain of a submanifold chart $(U _z,
\varphi)$ of
$P$. With this choice we can write
$U _z\simeq W _1 \times  W _2 $ and $V _z\simeq W _1\times\{0\} $,
where $W _1 $ and $W _2 $ are open neighborhoods of the origin in
two finite dimensional vector spaces of dimensions $\dim S $ and
${\rm codim}\, S $, respectively. If we denote the elements of 
$W_2$ by $(x _1, \ldots, x _{{\rm codim}\,S} )$ then any
arbitrary  extensions
$F _1, \ldots, F_{{\rm codim}\,S}\in  C^{\infty} (P) $ 
of the coordinate functions
$f_1= x _1 , \ldots, f_{{\rm codim}\,S}= x _{{\rm
codim}\,S}$ to the manifold $P$ generate
${\cal I}(V _z )$ and form a Poisson ideal. Indeed, since $V _z $
is an embedded quasi-Poisson submanifold of $P $, we
have for any $F \in  C^{\infty}(P) $ and any $z' \in V _z
$ 
\[
\{F _i, F\} (z')=\{f _i,F|_{V _z}\}_{V _z} (z')=0
\] 
since $f _i|_{V _z}\equiv 0 $.

Some of the ideas that we just introduced play a very important
role in the algebraic approach to Poisson geometry. The reader
interested in these kind of questions is encouraged to
check with~\cite{Vanhaecke} and references therein.

\begin{definition}
\label{I stability}
Let $(P, \{ \cdot , \cdot \}, H)$ be
a Poisson dynamical system and let $\mathcal{I}$ be a regular
Poisson ideal, that is, the vanishing set  $\mathcal{V_I}$ is a
closed and embedded submanifold of $P $.
Consider the reduced Poisson system $(\mathcal{V_I},
\{ \cdot , \cdot \}_{\mathcal{V_I}},h)$ where $h \in
C^{\infty}(\mathcal{V_I})$ is defined by 
$h:=H\circ i$ with $i:\mathcal{V_I}\hookrightarrow P$ the inclusion. 
Assume that $z_e\in
\mathcal{V_I}$ is an equilibrium point for the Poisson dynamical
system $(P, \{ \cdot , \cdot \},H)$ and hence also for
$(\mathcal{V_I},
\{ \cdot , \cdot \}_{\mathcal{V_I}},h)$. 
We say that $z_e\in \mathcal{V_I} \subset P$ is
an  $\mathcal{I}$-{\bfi stable  equilibrium} if any of
the two following equivalent conditions hold:
\begin{description}
\item[(i)] $z_e$ is a stable  equilibrium for the reduced
Poisson dynamical system $(\mathcal{V_I},
\{ \cdot , \cdot \}_{\mathcal{V_I}},h)$;
\item[(ii)] for any open neighborhood $U$ of $z_e$ in
$P$, there is an open neighborhood $V$ of $z_e$  in
$P $ such that if $F_t$ is the flow of $X_H$, then
$F_t(z)
\in U\cap \mathcal{V_I}$ for any $z \in V
\cap
\mathcal{V_I}$ and for any $t > 0$.
\end{description}
The equilibrium $z_e$ is ${\cal I}$-{\bfi  unstable}
if $z_e$ is an unstable  equilibrium for the reduced
Poisson dynamical system $(\mathcal{V_I},
\{ \cdot , \cdot \}_{\mathcal{V_I}},h)$. It is obvious that ${\cal
I}$-instability implies Lyapunov instability on the
whole space.
\end{definition} 

\begin{theorem}
\label{subcasimir stability} 
Let $(P, \{ \cdot , \cdot \},H)$ be a
Poisson dynamical system with an equilibrium at the point 
$z_e\in P$ and let $U \subset P$ be an open neighborhood
around $z_e$. Assume that there exists a regular Poisson
ideal $\mathcal{I}$ generated by the functions $G_1,
\ldots, G _m \in C^{\infty}(P)$ with sub-Casimirs $F_1,
\ldots , F_r \in C^{\infty}(P)$ such that $z_e \in
\mathcal{V_I}$. Suppose that the functions
$C_0:=H, C_1,\ldots,C_n \in C^{\infty}(P) $ are conserved
by the flow of $X _H $ and  that, additionally, there
exist constants $\lambda _1, \ldots, \lambda _n, \mu _1,
\ldots, \mu _r, \nu _1,
\ldots, \nu _m  $  such that 
\begin{description}
\item[(i)] $H_1:=\sum _{i=0}^{n}\lambda_i C_i  + \sum_{j=1}^{r}\mu_j
F_j+\sum_{k=1}^{m}\nu_k G_k$ has a critical point at
$z_e$, and
\item[(ii)] the Hessian of $H_1$ at $z_e$ is  positive
definite when restricted to  the subspace $W$ defined by
\[
W=
\bigcap_{i=0}^n\ker(\mathbf{d}C_i(z_e)) \bigcap_{j=1}^r
\ker(\mathbf{d}F_j(z_e))\bigcap_{k=1}^m
\ker(\mathbf{d}G_k(z_e)).
\]
\end{description} 
Then $z_e$ is an $\mathcal{I}$-stable equilibrium.
\end{theorem}
\textbf{Proof} The hypotheses in the statement of the theorem imply
that the equilibrium $z_e$ of the reduced system
$(\mathcal{V_I},\{ \cdot , \cdot \}_{\mathcal{V_I}},H
\circ i)$, with $i:\mathcal{V_I}\hookrightarrow P $  the
inclusion, satisfies the hypotheses of
Theorem~\ref{ultimate stability in Poisson} and hence is
Lyapunov stable on $\mathcal{V_I}$, which implies that
$z_e$ is ${\cal I}$-stable.
$\blacksquare$.

\section{Examples.}
\label{examples}
\subsection{A light Chaplygin sleigh on a cylinder.} 
The following
example was formulated in~\cite{nonholonomic} in the context of
nonholonomically constrained systems. In that work the author found
an equilibrium that exhibits asymptotically stable 
behavior. We will study the
stability of all the equilibria of this system as well
as of its relative equilibria with respect to a circle
symmetry of the system that will be introduced later on.
We will apply the Lyapunov stability methods presented
in the previous sections. This example is based on a 
real mechanical system that illustrates the theory
particularly well since it exhibits equilibria that are
not critical points of the Hamiltonian or of any  other
conserved function and, nevertheless,
Theorem~\ref{ultimate stability in Poisson} still allows
us to  establish the Lyapunov stability of some
dynamical elements and, in some cases, asymptotic
stability. There is also an equilibrium to which none of
the stability methods in the paper apply but 
that, after ideal reduction, is shown to be unstable
and hence unstable in the whole space.

We will start the presentation by explicitly carrying
out in this particular example Marle's reduction
procedure for nonholonomically constrained systems.  The
reader is encouraged to check with the original
references~\cite{nonholonomic,nonholonomic various} in order to
find various technical details that we will omit here.

\subsubsection{Description of the system.} 
The configuration space
is given by the points $(x,\theta) $ on a cylinder $
Q:=\mathbb{R}\times S^1$. The  Lagrangian of
the system is just the kinetic energy $L=\frac{1}{2}(\dot{x}^2 +
\dot{\theta}^2)\in C^{\infty}(TQ)$. The system
is constrained to move subject to the semiholonomic
constraint $\dot{x}+x \dot{\theta}=0$. The term
``semiholonomic" means that the distribution that
describes the constraint is integrable with integral
leaves that are not necessarily embedded submanifolds.

This system  approximates  a simple
mechanical system in a certain regime that can be physically
realized in the following way. Take a Chaplygin sleigh moving in
the interior of a cylinder (we are assuming that all the physical
constants of the system are equal to $1$). The configuration space
of this system consists of the points 
$(x,\theta,\varphi) \in Q':=\mathbb{R} \times S^1 \times S^1$, where
the coordinates $(x,\theta)$  on the cylinder indicate the
position of the Chaplygin sleigh.  The dynamics of this
system is determined by the Lagrangian $L'  $  on $TQ' $ given by
$L'=\frac{1}{2}(\dot{x}^2 +
\dot{\theta}^2+I_{\varphi}\dot{\varphi}^2)$, where $I_{\varphi}$ is
the moment of inertia of the sleigh, together with the nonholonomic
constraint $\dot{x}\cos \varphi -\dot{\theta} \sin{\varphi}=0$.
Assume now that we add a new holonomic constraint $\tan
\varphi=x$. Notice that even if the first constraint was not
integrable, the superposition of the two constraints 
\textit{is}  integrable. In this case the  dynamics
can be described  by restricting the system to a new
configuration space $\bar{Q}
\subset Q'$ which is actually an integral manifold of the
distribution that describes the holonomic constraint.
Moreover, it is easy to see that we can restrict the
system to the integral manifold of
\textit{any} subset of integrable constraints, obtaining a new
holonomically constrained system. In this case, we restrict the
system described by the Lagrangian $L' $ on $TQ' $  to the
integral submanifold
$Q\subset Q'$ by using the holonomic constraint $\tan
\varphi=x$. Assuming $I_{\varphi}\ll 1$ and restricting our study
to points such that  $x \ll 1$, the example that  we will be
presenting is a good approximation of this mechanical system.
Marle~\cite{nonholonomic} considers the same mechanical realization
of these equations  but he sets  $I_{\varphi}=0$ from the
beginning of his exposition. 

\subsubsection{Reduction of the system} 
We now apply  a reduction
procedure  due to Marle~\cite{nonholonomic, nonholonomic
various} to eliminate the semiholonomic constraint $\dot{x}+x
\dot{\theta}=0$. This reduction procedure consists of eliminating
the Lagrange multipliers of a (in general nonholonomically)
constrained system by finding a submanifold (the constraint
submanifold) endowed with an almost Poisson structure and a
Hamiltonian on it in such a way that the dynamics of this almost
Poisson dynamical system coincides with the dynamics of the
original constrained system. There are several equivalent
constructions (see~\cite{van der Schaft and Maschke 1994,cushman
kemppainen sniatycki,nonholonomic,nonholonomic mechanical with
bloch,sniatycki 2001}, and references therein) to handle these
constraints. It was shown in~\cite{van der Schaft and Maschke
1994} that this almost Poisson structure is actually Poisson if and
only if the constraints are semiholonomic.

Let $Q=\mathbb{R} \times S^1$ be the configuration space
and  $L(x,\theta,
\dot{x},\dot{\theta})=\frac{1}{2}(\dot{x}^2 +
\dot{\theta}^2)$ the Lagrangian of the system subjected
to the constraint $\dot{x}+x\dot{\theta}=0$. Since the
Lagrangian $L$ is hyperregular, the Legendre transform  
$\mathbb{F}L:TQ \rightarrow T^{\ast}Q$ is an
isomorphism that we use to associate  a Hamiltonian
function
$H \in C^{\infty}(T^{\ast}Q)$ to the system. The image by
$\mathbb{F}L $ of the constraint submanifold in
$TQ$ gives the constraint submanifold $P$ on
$T^{\ast}Q$ which consists of the points
$P=\{(x,\theta,p_x,p_{\theta})
\in T^{\ast}Q \mid  p_x+x p_{\theta}=0 \}$.
Let $\mathcal{D}\subset T(T^{\ast}Q)$ be the so called {\bfi
constraint distribution}  defined by
$\mathcal{D}(z):=T_zP$ for every $z$ in $P$.
D'Alembert's principle provides a prescription to modify
the  original unconstrained Hamiltonian flow in order to
construct a new vector field whose integral curves lie
in $P$. Indeed, let
$\left. X_H\right|_{P}$ be the restriction of the original
Hamiltonian flow to the points in $P$ and let $X_D$ be the modified
vector field whose integral curves describe the dynamics of the 
nonholonomically constrained system. The works by Marle
quoted above  ensure that, under certain regularity
conditions satisfied in this example, the difference
$X_W=\left.X_H\right|_{P}-X_D$ of these two vector
fields, is a section of  a  subbundle ${\cal W}$ of
$T_{P}(T ^\ast Q)$ that satisfies
$T_P(T ^\ast Q)={\cal W}\oplus \mathcal{D}$ and that is
uniquely determined  by D'Alembert's principle. In such a
situation, every Hamiltonian vector field can be
decomposed in a unique way as
$\left.X_H\right|_{P}=X_D+X_W$ and $X_D$ describes the dynamics of
the constrained system. Marle also shows that there
exists an almost Poisson structure on $P$ with almost Poisson
tensor
$B:T^{\ast}P
\times T^{\ast}P \ \longrightarrow \mathbb{R}$,  for which
$X_D=B^{\sharp} \mathbf{d}H|_{P}$, where $B^{\sharp}:T ^\ast P
\rightarrow TP$ is the canonical vector bundle isomorphism
associated to $B$.

In our example, $\mathcal{D}(x,\theta,p_x,p_{\theta})={\rm
span}\{(1,0,-p_{\theta},0),(0,1,0,0),(0,0,-x,1)\}$ and 
$\mathcal{W}(x,\theta,p_x,p_{\theta})={\rm span} \{(0,0,1,x) \}$.
An explicit expression for the almost Poisson
structure (see~\cite{ortega planas 2004}) can be given by using the
natural  projection map  onto the $\mathcal{D}$ factor. After some
computations this almost Poisson tensor takes the form:
\[ B(x,\theta,p_{\theta})=\begin{pmatrix} 0 & 0 & \frac{-x}{1 +
x^2}\\ 0 & 0 & \frac{1}{1+x^2}\\
\frac{x}{1+x^2}&\frac{-1}{1+x^2}&0
\end{pmatrix}
\] where the three-tuples $(x,\theta,p_\theta)$ are used to
coordinatize the points $(x,\theta,-x p _\theta, p_\theta)\in P$ and
the restricted Hamiltonian is given by
$H|_P(x,\theta,p_{\theta})=\frac{1}{2}(1+x^2)p_{\theta}^2$. Notice
that this tensor is Poisson since the constraint is integrable.
The equations of motion are
\[
\dot{x}=-x p_{\theta} \quad,\quad 
\dot{\theta}=p_{\theta}\quad,\quad 
\dot{p}_{\theta}=\frac{x^2 p_{\theta}^2}{1+x^2}\,.
\]

\subsubsection{Equilibria, relative equilibria, and their stability}

Notice that every point of the form $z=(x,\theta,0)$ is
an equilibrium of the system. If we first compute the
linearization of the dynamical system at those equilibria we obtain
the family of matrices
\[ 
\begin{pmatrix} 0 & 0 & -x\\ 0 & 0 & 1\\ 0 & 0 & 0
\end{pmatrix}
\] 
which have three zero eigenvalues and are not diagonalizable.
This implies that the system   is linearly unstable at those
equilibria (which does not imply either Lyapunov
stability or instability). 

To apply Theorem~\ref{ultimate stability in Poisson}, we
first need to find conserved quantities for the
Hamiltonian flow.  In this case we can use
the Hamiltonian and the local Casimir function given by
$C(x,\theta,p_\theta)=x e^{\theta}$. Let $L$ be the function
defined by
$L :=\lambda_0 H +
\lambda_1 C$. If we set $\lambda_0=1$ and $\lambda_1=0$ we have
that  $\mathbf{d}L(z)=0$. The subspace
$W=\ker
\mathbf{d}H(z) \cap \mathbf{d}C(z)$ is given by $W={\rm
span}
\{(x,-1,0),(0,0,1) \}$ and the restricted Hessian
\[
\left.\mathbf{d}^2L(z)\right|_{W\times W} =
\begin{pmatrix} 0 & 0\\ 0 & (1 + x^2)
\end{pmatrix}
\]
is not positive definite since it has a zero eigenvalue.
The stability of the equilibrium $z=(0,\theta,0)$ can be
analyzed by using the fact that the submanifold $S$
consisting of the points of the form
$(0,\theta,p_{\theta})$ is such that  its vanishing ideal
$\mathcal{I}(S)  $ is a Poisson ideal and hence $S$ is Poisson
reducible. Indeed, if $(\theta,p_{\theta})$ are coordinates on 
$S$, the reduced bracket $\{ \cdot , \cdot \} _S $ takes
the form $\{\theta,p_{\theta}\} _S=1$ and the reduced
Hamiltonian is
$h(\theta,p_{\theta})=\frac{1}{2}p_{\theta}^2$. This
reduced system describes a free one dimensional
particle. The equilibrium $z=(0,\theta,0)$ of the
original system drops to an equilibrium at the point
$(\theta,0)$ which is clearly unstable. In particular,
this implies the instability of the original equilibrium
$(0,\theta,0)$.

We now study the stability of the relative equilibria with respect
to  the circle symmetry of the system given by the 
action $\psi \cdot
(x,\theta,p_{\theta})=(x,\theta+\psi,p_{\theta})$. This action is
canonical and the system can be Poisson reduced. The reduced
manifold is $\mathbb{R}^2$. If we denote by
$(x,p_{\theta})$ the elements of the reduced space, the
reduced Poisson bracket is determined by the relation
$\{x,p_{\theta}\} =-x/(1+x^2)$ and the reduced
Hamiltonian is
$h(x,p_{\theta})=\frac{1}{2}(1+x^2)p_{\theta}^2$.
Hamilton's equations for $h $ are 
$\dot{x}=-x p_{\theta},\;
\dot{p}_{\theta}=x^2 p_{\theta}/(1+x^2)$. Thus the
equilibria are given by the family of points satisfying
$xp_{\theta}=0$. The linearization of the Hamiltonian
vector field at these equilibria is given by the matrix
\[ 
\begin{pmatrix} -p_{\theta}& -x\\ 0 & 0
\end{pmatrix}
\]
which has a positive eigenvalue if $p_{\theta}<0$, in
which case the system is Lyapunov unstable at the points
$(0,p_{\theta})$. This obviously implies that the
unreduced system exhibits nonlinearly unstable relative
equilibria.

If $p_{\theta}>0$ the linearization does not imply 
neither stability nor instability. However, note that in
this case, the linearization has a negative eigenvalue
with eigenvector $v=(1,0)$ that will be useful when
searching for a Lyapunov function (see
Remark~\ref{lyapunov function and negative
eigenvalues}). In order to study the nonlinear stability
of these relative equilibria, we notice that the only
available conserved quantity is the reduced Hamiltonian
whose derivative
$\mathbf{d}h(x,p_{\theta})=(x p_{\theta}^2,
(1+x^2)p_{\theta}) = (0,0)$ 
if and only if
$p_{\theta}=0$. In that case 
\[
\mathbf{d}^2H(x,0)= \begin{pmatrix} 0 & 0\\ 0&
(1+x^2)
\end{pmatrix} 
\] 
and hence  we cannot conclude
either stability or instability. However, in this
particular case instability can be concluded just by
looking at the phase portrait for the vector field.
For points of the form $(0,p_{\theta})$ the derivative
of the Hamiltonian does not vanish and hence the only
way to apply Theorem~\ref{ultimate stability in Poisson}
consists of finding a function $F$
satisfying at least one of the hypotheses {\bf (i)} or
{\bf (ii)};  
$F(x,p_{\theta})=x^2/2$ is one such function since $\{x^2,h\}=-2x^2
p_{\theta}$,
$\{x^4,h\}=-4x^4 p_{\theta}$, and
$p_{\theta}$ is assumed to be positive. Consequently, the hypothesis
{\bf (i)} is obviously satisfied. With this choice, the subspace
$W=\ker
\mathbf{d}h(0,p_{\theta}) ={\rm span}\{(1,0) \}$ and
$\mathbf{d}^2 F(0,p_{\theta})|_{W \times W}=1>0$. Consequently, the
equilibria of the form
$(0,p _\theta)$ with $p_\theta>0$ are 
Lyapunov stable and even though they are not asymptotically
stable, there exists an open neighborhood $V$ of
$(0,p_{\theta})$  such that $F _t(V)\subset F _s(V) $, whenever $t>s
$, that is, they are weakly asymptotically stable.

Finally, it is easy to conclude that the equilibria on 
the form $(x,0)$ are unstable just by looking at the
phase portrait of the system.

\subsection{Two coupled spinning wheels}
\label{Two coupled spinning wheels}

Consider two vertical weightless wheels with radii $R$ 
and $r$ satisfying $R> r$ and 
$R/r\in \mathbb{R}\setminus \mathbb{Q}$. We attach to
the edges of  each of these wheels two point masses
$M$ and $m$. This simple system has as
configuration space $Q$ the torus $\mathbb{T}^2$ that we
coordinatize with the angles $(\theta, \varphi)$. The
Lagrangian of this system in these coordinates is
$L=\frac{1}{2}(M R^2 \dot{\theta}^2 + m
r^2 \dot{\varphi}^2)+ M R \cos \theta + m r \cos \varphi$.
Assume now that we couple the rotations of the two wheels
with a belt. This mechanism imposes on the systems a
semiholonomic constraint that can be expressed as $R
\dot{\theta} - r \dot{\varphi}=0$. In order to give a
description of the constrained system  we first express
the original system in the Hamiltonian setting by using
the Legendre transform. The phase space $P$  is in this
case the cotangent bundle
$T^{\ast}\mathbb{T}^2 \simeq \mathbb{T}^2 \times \mathbb{R}^2$
with coordinates $(\theta, \varphi,p_{\theta},p_{\varphi})$,
endowed with the canonical symplectic form. The Hamiltonian
function is 
\[H=\frac{1}{2}\left(\frac{p^2_{\theta}}{M R^2}+
\frac{p^2_{\varphi}}{m r^2}  \right) - M R \cos \theta - m r \cos
\varphi.\] The constraint submanifold is given by the
points$(\theta, \varphi, p_{\theta}, p_{\varphi}) $ that
satisfy $p_{\varphi}= m r p_{\theta}/MR$, which can be
identified with $\mathbb{T}^2 \times \mathbb{R}$ with
coordinates $(\theta, \varphi, p)$.
\begin{figure}[htb]
\begin{center}
\includegraphics{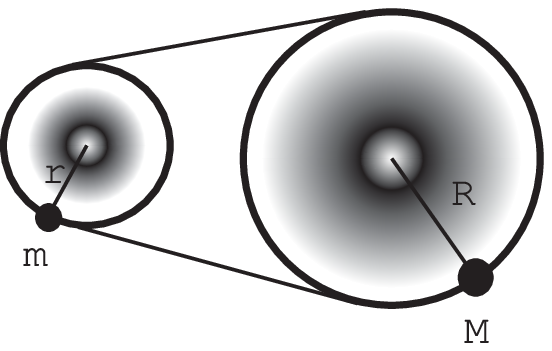}
\end{center}
\end{figure}

We now apply the  reduction procedure in~\cite{Bates} in order
to find a bracket on the constraint submanifold that 
is actually Poisson since the constraint is
semiholonomic. This bracket is given by the
constant Poisson tensor:
\[
B(\theta, \varphi,p)=\begin{pmatrix}
0 & 0 & r\\
0 & 0 & R\\
-r & -R & 0
\end{pmatrix}.
\]
The reduced Hamiltonian function is 
\[h(\theta,
\varphi,p)=\frac{p^2}{2 k}-M R \cos
\theta- m r
\cos  \varphi,\] 
where $k$ is a real positive constant depending on
the parameters of the problem given by the expression 
\[k=\frac{m +
M}{4 M m R^2 r^2}-\frac{(m-M)^2m^2 M^2}{4 m^2 M^2 R^2 r^2 (m+M)}.\]
This Poisson system has a local Casimir given by the locally defined
function
$C(\theta, \varphi, p)=R \theta - r \varphi$.
The equations of  motion of the system are given by
\[
\dot{\theta}=r \frac{p}{k}, \qquad \dot{\varphi}=R \frac{p}{k},
\qquad \dot{p}=-rRM \sin \theta -m r R \sin \varphi.
\]
The equilibria of the system are the points of the set
$S=\{(\theta,
\varphi, 0)\mid M \sin \theta + m \sin
\varphi=0\}$ that can be described as a one-parameter family
given by the curve
$\varphi=-\sin^{-1}\left(\frac{M \sin \theta}{m}\right), \; \theta
\in [-\theta_c , \theta_c]$, where $\theta_c$ is given by
$\theta_c=\sin^{-1}(\frac{m}{M})$. In order to study the nonlinear
stability of such equilibria we compute $ \lambda_0\mathbf{d}h(z) +
\lambda_1 \mathbf{d}C(z)=0$, with $z=(\theta,
\varphi, 0) \in S $. This equation can be
solved by taking
$\lambda_1=-\lambda_0 M \sin \theta$. In this case $W=\ker
\mathbf{d}C(z)
\cap
\ker \mathbf{d}h(z)={\rm span} \{(r,R, 0),(0,0,1)\}$. Finally, it
is easy to see that 
$\mathbf{d}^2(\lambda_0h + \lambda_1 C)(z)|_{W \times W}>0$ if and
only if $Mr \cos \theta + m R \cos \varphi>0$. In particular, the
point $z=(0,0,0)$ is always nonlinearly stable, as
expected, and the point $z=(0, \pi, 0)$ is stable if
$\frac{M}{R}> \frac{m}{r}$.

\section{Nonlinear stability via topological methods}
\label{how to compare withg Patrick}
\normalfont

In \cite{topological stability} topology based tools
have been developed that provide sufficient conditions
for the Lyapunov stability of Poisson equilibria. One of
the main achievements in~\cite{topological
stability}  is the discovery of a space related
to the topology of the symplectic foliation of the Poisson manifold
(see~(\ref{t2 over line}) below) on which the extremality of the
Hamiltonian suffices to conclude stability. 
In this section we will study under which circumstances the
topological criteria in~\cite{topological
stability} can be expressed in terms of local continuous Casimir
functions and hence there is an equivalence with the
energy-Casimir method. To be more explicit, we will seek
the correspondence between the topological approach of
\cite{topological stability} and a generalization of the
energy-Casimir method that requires only continuity of
the functions involved and that is based on the
following general lemma.

\begin{lemma} Let $X \in
\mathfrak{X}(P)$ be a smooth vector field on the finite dimensional
manifold $P$ and $z_e \in P$ an equilibrium point. If
there exists locally defined continuous conserved
quantities $C_0, \ldots , C_k
\in C^{0}(U)$ of the flow $F _t $ such that $\bigcap_{i=0}^k
C_i^{-1}(C _i (z_e))=\{z_e\}$  then the equilibrium $z_e $
is Lyapunov stable.
\end{lemma}
\noindent
\textbf{Proof} Consider the function
$L(z)= (C_0(z)-C _0(z_e)) ^2 + \ldots +
(C_k(z)-C _k (z_e)) ^2 $. The hypothesis $\bigcap_{i=0}^k
C_i^{-1}(C _i (z_e))=\{z_e\}$ ensures that $L$ is a
positive function that takes the zero value only at the
point $z_e $. In particular, the sets of the form $L
^{-1}([0, \epsilon)) $,
$\epsilon>0 $, form a fundamental system of
neighborhoods in the manifold topology of
$P$ at the point $z_e$. Consequently, for any open
neighborhood $U $ of  $z_e$ there exists an $\epsilon>0 $
such that  $L ^{-1}([0, \epsilon])
\subset U$. Since the level set $L ^{-1} (\epsilon) $
is invariant by the flow $F _t  $ of $X$, the Lyapunov
stability of $z_e$ follows.
$\blacksquare$
\medskip

Any continuous function $C \in  C^0 (U) $, with $U $ an
open subset of $P$, such that  $C$  is constant on the
symplectic leaves of $(U,\{ \cdot ,\cdot\}| _U )$
is called a  {\bfi  local continuous Casimir} of $(P,\{
\cdot , \cdot \})$. The choice of terminology is
justified by the fact that if such a function $C $ happens
to be differentiable then it is an actual Casimir of
$(U,\{ \cdot ,\cdot\}| _U )$. It is worth noticing that
the local continuous Casimirs are the (continuous) first
integrals of the foliation of $(U,\{ \cdot ,\cdot\}| _U
)$ by its symplectic leaves.

\begin{corollary}[Continuous energy-Casimir method] 
\label{how to compare both}
Let $(P,\{ \cdot , \cdot \},H)$ be a Poisson
dynamical system and $z_e \in P$  an equilibrium point
of the Hamiltonian vector field $X_H$. Let $S_{z_e}
\subset P$ be the common level set of  local continuous
Casimir functions around $z_e$. If
\begin{equation}
\label{continuous energy casimir}
H^{-1}(H (z_e))\cap S _{z_e}=\{z_e\}
\end{equation}
then the equilibrium $z_e$ is Lyapunov
stable. This statement remains true if $H$ is replaced by any
continuous conserved quantity of the flow of $X _H$.
\end{corollary}

Our goal is to establish sufficient conditions under which this
corollary coincides with the topological stability criterion
in~\cite{topological
stability} that we now recall. We start by introducing the
necessary notation. Let~$(X,\tau)$ be a topological space and $x\in X$
an arbitrary  point. We define the set $T_2(x)\subset X$ as
\begin{equation}
\label{definition t2 sets} 
T_2(x):=\{y \in X \mid  U_x \cap U_y
\not=\emptyset \text{ for any two open neighborhoods }U_x, \, U_y
\text{ of  $x$ and $y$}\}.
\end{equation}
Let $A \subset X$ be an arbitrary subset. We define
\[
T_2(A):=\bigcup _{x \in A} T _2 (x).
\]
Notice that if $y \in T_2(x)$ then $x \in T_2(y)$. Also, a
topological space
$(X,\tau)$ is Hausdorff if and only if $T_2(x)=x$ for every $ x \in
X$. Hence the $T _2 $ sets measure how far a topological space is
from being Hausdorff.

Suppose now that $P$ is a smooth Hausdorff and paracompact finite
dimensional manifold and $D$ is a smooth and integrable generalized
distribution on $P$. Let $\pi_D:P \rightarrow P/ D $ be the
projection onto the leaf space of the distribution $D$. The map
$\pi$ is continuous and open when $P/ D  $ is endowed with the
quotient topology. Define
\begin{equation}
\label{t2 over line}
\overline{T}_2(x)=\pi_{D}^{-1}\left(
T_2\left(\pi_{D}(x)\right) \right),\quad x\in P,
\end{equation}
and, more generally, 
\begin{equation}
\label{t2 over line with U}
\overline{T}_2^U(x)=\pi_{D| _U}^{-1}\left(
T_2\left(\pi_{D| _U}(x)\right) \right),\quad x\in P,
\end{equation}
where $U $  is an open neighborhood of $x \in  P $ and $\pi_{D|
_U}:U \rightarrow  U/ D| _U $ is the projection onto the leaf space
of the restriction $D | _U  $ of $D$ to $U$.

We now focus on the particular case when $P $ is a Poisson
manifold with bracket $\{ \cdot , \cdot  \}$. Let $ \mathcal{E} $
be the corresponding characteristic distribution and
$\pi:P \rightarrow P/ \{ \cdot , \cdot \}$ the projection
onto the space of symplectic leaves $P/ \{ \cdot , \cdot \}:=P/
\mathcal{E} $.

\begin{theorem}[Topological energy-Casimir method~\cite{topological
stability}] 
\label{patrick theorem statement}
Let $(P,\{ \cdot , \cdot \},H)$ be a Poisson
dynamical system and $z_e \in P$  an equilibrium point
for the Hamiltonian vector field $X_H$. If there is an
open neighborhood $U \subset P$ of  $z_e$ such that  
\begin{equation}
\label{t 2 energy topological}
H^{-1}(H (z_e))\cap \overline{T}_2^U(z_e)=\{z_e\}
\end{equation}
then the equilibrium $z_e$ is Lyapunov
stable. This statement remains true if $H$ is replaced by any
continuous conserved quantity of the flow of $X _H$ that takes
values in a Hausdorff space.
\end{theorem}

In view of expressions~(\ref{continuous energy casimir})
and~(\ref{t 2 energy topological}) we would like to know under what
circumstances the set $ \overline{T}_2^U(z_e) $ can be
obtained by looking at the level sets of local
continuous Casimir functions thereby rendering the
statements of Corollary~\ref{how to compare both} and
Theorem~\ref{patrick theorem statement} equivalent. 

\begin{figure}[htb]
\begin{center}
\includegraphics{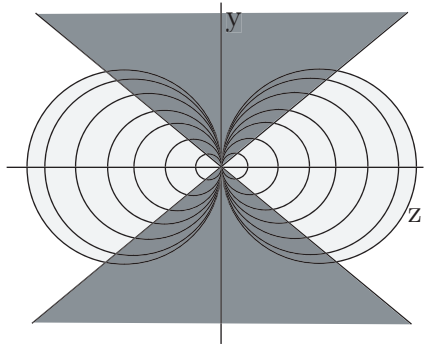}
\end{center}
\caption{Symplectic leaves of Montaldi's example of a Poisson
manifold that does not have local Casimirs around the origin. The
shadowed area represents the set $\overline{T} ^U_2(0,0,0) $. The
picture is a section of the three dimensional figure through the
$OYZ$ plane.}\label{fig:james example}
\end{figure}

The first point
that we have to emphasize is that this is, in general, not
possible. The following example, that we owe to James Montaldi, shows
that, in general, we cannot find enough local Casimir functions
to be able to write the set $ \overline{T}_2^U(z_e) $ as the common
level set of local continuous Casimir functions,
no matter how much we shrink the neighborhood $U$. Let $\mathbb{R}^3$ 
and $f(x,y,z) = x^2 + y^2 - z^2$. Consider the Poisson
structure $\{ \cdot , \cdot \} $ determined by 
$\{x, y\} = f^2 $,   $\{y,z\} =  2yzf $, and  $ \{x,z\} = -2xzf$.
In order to  describe the symplectic leaves of $(\mathbb{R} ^3,\{ \cdot
, \cdot \} )$ (see Figure~\ref{fig:james example}) notice first that
the function $f $ is a factor and hence the Poisson tensor vanishes on
the cone $f=0$. Consider now
all the spheres through the origin and tangent to the $OXY $ plane
(and hence centered on the $OZ $-axis) and  cut them with the
cone $f = 0$. Each of these spheres contains the
following symplectic leaves: the sphere intersected with
the points $(x,y,z) $ such that 
$f(x,y,z)>0 $ (two dimensional leaf), the sphere intersected with 
the points $(x,y,z) $ such that 
$f(x,y,z)<0 $ (two dimensional leaf), and the points such that 
$f(x,y,z)=0$ (zero dimensional leaves).
It is clear from this description that there are no non-constant
continuous local Casimir functions near the origin. Nevertheless,
for any neighborhood $U $ of the origin
$\overline{T} ^U_2(0,0,0) = \{(x,y,z)\mid f(x,y,z)  \geq 0\} $,
that is, the closed exterior of the cone, which in this case is
strictly included in $C _U ^{-1}(C _U(0,0,0))= U $.

Even though the previous example shows that the set 
$ \overline{T}_2^U(z_e) $ does not coincide in general with the
common level set of local continuous Casimir functions one can easily
prove that at least  one inclusion holds true. The natural context to 
present most of the results in this section is that of
generalized foliations of smooth manifolds. Consequently, we
will prove our statements in that category and we will
obtain the Poisson case as a corollary by applying the theorems
to the generalized foliation of the Poisson manifold by its
symplectic leaves.

\begin{lemma}
\label{lemma famous inclusion for ever}
Let $P$ be a smooth finite
dimensional manifold and $D$  a smooth integrable
generalized distribution on $P$. Let $\pi_D:P
\rightarrow P/ D $ be the projection onto the leaf space
of the distribution $D$ and
$\overline{T}_2 $ the  symbol  defined in \eqref{t2 over
line}. Let $C _i \in  C ^0 (P) $, $ i \in  I $, be a set of continuous functions
that
are constant on the integral leaves of $D$ (that is, 
first integrals of
$D$). Then for any $z \in P $
\begin{equation}
\label{famous inclusion for ever}
\overline{T}_2 (z )\subset \bigcap _{i\in  I} C _i
^{-1}(C _i (z )).
\end{equation}
\end{lemma}

\noindent\textbf{Proof.\ \ } Let
$C: P
\rightarrow
\Bbb R ^I$ be the  function defined by $C (z):=(C _i (z))_{i \in  I}
$. If we endow $\mathbb{R}^I $ with the product topology (not the box
topology!) then the continuity of the first integrals $C _i$, $ i \in
I $, implies that $C$ is continuous. The projection
$\pi:P \rightarrow  P/D$ is an open map when $P/D$ is endowed with the
quotient topology. Given that $C$ is constant on the integral 
leaves of  $D$ it drops to a map $c :  P/ D \rightarrow \Bbb R ^I$
that closes the diagram
\unitlength=5mm
\begin{center}
\begin{picture}(9,6)
\put(1,5){\makebox(0,0){$P$}} 
\put(9,5){\makebox(0,0){$\Bbb R^I$}}
\put(5 ,0){\makebox(0,0){$P/D$}} 
\put(1,4){\vector(1,-1){3}}
\put(6,1){\vector(1,1){3}} 
\put(2,5){\vector(1,0){6}}
\put(1,2){\makebox(0,0){$\pi_D $}}
\put(9,2){\makebox(0,0){$c $}}
\put(5,6){\makebox(0,0){$C $}}
\end{picture}
\end{center}
The continuity of $C $ and the openness and 
surjectivity of $\pi_D $
imply that
$c $ is also continuous.  In order to
prove~(\ref{famous inclusion for ever}) it suffices to
show that  if  $m \in \overline{T} _2   (z) $ then $C
  (m)= C   (z )$. By contradiction, suppose that $ C
 (m)\neq C   (z )$. Since $ \Bbb R^I$ is a
Hausdorff topological space there are open neighborhoods
$V _{C (m)}$ and $V _{C  (z )}$ of $m$ and  $z $,
respectively, such that $V _{C  (m)}\cap V _{C(z )}=
\emptyset$. As $c  $ is continuous the sets 
$c^{-1}(V_{C(m)}) $ and $c^{-1}(V _{C(z)})$ are open
neighborhoods of $\pi _D (m)$ and $\pi_D
 (z )$, respectively. Also, since by hypothesis
$m \in \overline{T} _2   (z )$, we have that $c 
^{-1}(V _{C   (m)}) \cap c   ^{-1}(V _{C 
(z )}) \neq \emptyset$, necessarily. However, by
construction we also have that $c   ^{-1}(V _{C  
(m)}) \cap c   ^{-1}(V _{C  
(z )})= c   ^{-1}(V _{C  
(m)}\cap V _{C  
(z )})= c  ^{-1}(\emptyset)=\emptyset$, which is a
contradiction. \quad $\blacksquare$

\medskip

The rest of this section is dedicated to the description of
two situations where the inclusion~(\ref{famous inclusion for
ever}) is an equality and hence local continuous Casimir functions
characterize the $\overline{T} _2 $-sets. 
 We
start with a couple of preliminary general results.

\begin{definition}
Let $(X, \tau)$ be a topological space. We say that $(X, \tau)$ is
$T _2 $-{\bfi  idempotent} when $T _2(T _2 (x))= T _2 (x) $, for
any $x \in X $.
\end{definition} 

\begin{lemma}
\label{equivalence relation} 
Let $(X,\tau)$ be a $T _2
$-idempotent topological space.
Then 
\begin{description}
\item[(i)] The relation $ \mathcal{R}_{T_2} $ on $X$ defined by
$x\mathcal{R}_{T_2} y$ if and only if $y \in T_2(x)$ is
an equivalence relation on $X$.
\item[(ii)]  The following statements are equivalent:
\begin{enumerate}
\item  $y \notin T_2(x)$.
\item  $T_2(x) \not= T_2(y)$.
\item  $T_2(x) \cap T_2(y)=\emptyset$.
\item  There exist  open neighborhoods $U_x,\,U_y$ of $x$ and
$y$, respectively, such that $T_2(U_x) \cap T_2(U_y)=\emptyset$.
\end{enumerate}
\item[(iii)]  If the projection $\pi_{T_2}:X \longrightarrow
X/\mathcal{R}_{T_2}$ onto the space of equivalence classes endowed
with the quotient topology is an open map then 
$X/\mathcal{R}_{T_2}$  is a
Hausdorff topological space.

\end{description}
\end{lemma}

\noindent
\textbf{Proof} {\bf (i)} The definition of the $T _2 $ 
set implies that $x \mathcal{R}_{T_2} x$ for any $x \in
X $ and that $x \mathcal{R}_{T_2}y$ if and only if $y
\mathcal{R}_{T_2} x$. In order to prove transitivity
of $\mathcal{R}_{T_2} $ let $x,y,z \in X $ be such that 
$x\mathcal{R}_{T_2} y$ and $y \mathcal{R}_{T_2} z$. 
By the very definition of the $T _2 $ set, it is clear 
that for any two subsets $A, B \subset X$ such that  $A
\subset B $ we have that $T _2(A) \subset  T_2(B)$, in
particular, the condition
$x
\in T_2(y)$ implies that
$T_2(x)
\subset T_2 (T_2(y))=T_2(y)$. By reflexivity we have that $T_2(y)
\subset T_2(x)$ and hence $T_2(x)=T_2(y) $ which implies that
$T_2(x)=T_2(y)=T_2(z)$ and hence $x \mathcal{R}_{T_2} z$.  

\smallskip

\noindent {\bf (ii)} If $T_2(x)=T_2(y)$ then $y \in T_2(y)=T_2(x)$.
This proves the implication
1$\Rightarrow$2. The implication 2$\Rightarrow$1 was
already proved in the first part of the lemma.  In
order to prove 2$\Rightarrow$3 suppose that there exists a point $z
\in T_2(x)
\cap T_2(y)$. Then using the $T _2 $ idempotency as we
did in the proof of the first part of the lemma we
obtain that
$T_2(x)=T_2(z)=T_2(y)$, which contradicts the hypothesis.
To show 3$\Rightarrow$4, assume that
$T_2(x) \cap T_2(y) = \emptyset$. Then, in particular,
$y\notin T _2 (x) $ and hence there exist  open
neighborhoods $U_x$ and $U_y$ of
$x$ and
$y$, respectively, such that $U_x \cap U_y =\emptyset$. Since
$U_x$ and $U_y$  are open neighborhoods of each of
their points, it follows that for every
$a_x\in U_x$ and
$a_y\in U_y$ the element $a _x \notin T _2(a _y )$. Using the
implication 1$\Rightarrow$3 that we have already proved, this
shows  that $T _2(a _x)\cap T _2(a _y )=\emptyset $ and hence 
\[
T_2(U_x) \cap T_2(U_y) =\left( \bigcup _{a_x \in U _x}T
_2(a _x)\right)\bigcap \left( \bigcup _{a_y \in U _y}T
_2(a _y)\right)=\bigcup_{a_x \in U _x,a_y \in U
_y}\left(T _2(a _x) \bigcap  T _2(a
_y)\right)=\emptyset.
\]
Finally, the implication 4$\Rightarrow$2 is straightforward.

\smallskip

\noindent {\bf (iii)} Notice first that for every subset $A \subset
X$, we have that
$\pi^{-1}_{T_2}( \pi_{T_2}( A))=T_2(A)$.  Let $\rho, \sigma
\in X/\mathcal{R}_{T_2}$ be two  points such that  $\rho\neq \sigma
$ and let
$x$ and
$y$ be two points in $X$ such that $\pi_{T_2}(x)=\rho$ and
$\pi(y)=\sigma$. Since $T_2(x) \not= T_2(y)$  there
exist, by part {\bf (ii)}, two  open neighborhoods $V _x $  and $V
_y 
$ of 
$x $ and 
$y
$, respectively, such that   
$\emptyset=T_2(V_x) \cap T_2(V_y)=\pi_{T_2}^{-1}(\pi_{T_2}( V_x))
\cap \pi_{T_2}^{-1}(\pi_{T_2}( V_y))=\pi_{T_2}^{-1}(\pi_{T_2}( V_x)
\cap \pi_{T_2}( V_y))$. Applying $\pi _{T _2  } $  to both sides of
this equality we obtain that
$\pi_{T_2} (V_x) \cap \pi _{T _2}( V_y) =\emptyset$. Since
$\pi_{T_2} (V_x)$ and
$\pi_{T_2} (V_y)$ are, by the openness of $\pi _{T _2} $, open
neighborhoods of the points $\rho$ and $\sigma$,
respectively, the claim follows.
$\blacksquare$

\medskip

Suppose now that $P$ is a smooth Hausdorff and paracompact finite
dimensional manifold and $D$ is a smooth and integrable generalized
distribution on $P$. Let $\pi_D:P \rightarrow P/ D $ be the
projection onto the leaf space of the distribution $D$ and
$\overline{T}_2 $ the  symbol defined in~(\ref{t2
over line}). Notice that since $\pi _D $ is surjective,
we have  
\begin{equation}
\label{pi equivariance t 2 e}
\pi_D(\overline{T}_2(x))=
T_2\left(\pi_{D}(x)\right),\text{ for any } x\in P.
\end{equation}
We will say that the pair $(P, D)$ is $\overline{T}_2 ${\bfi 
-idempotent} when $\overline{T}_2(\overline{T}_2
(x))=\overline{T}_2 (x) $, for any $x \in  P $. Notice that since the
sets
$\overline{T}_2(x)
$ are
$D$-saturated (they are unions of leaves of $D$), we can
conclude using~(\ref{pi equivariance t 2 e}) that $P$ is
$\overline{T}_2$-idempotent if and only if $P/ D $ is $T _2
$-idempotent. With this remark in mind the previous lemma can
be easily adapted to the symbol $\overline{T}_2$.

\begin{lemma}
\label{lemma upstairs}
Let $P$ be a smooth Hausdorff paracompact finite
dimensional manifold and $D$  a smooth integrable
generalized distribution on $P$. Let $\pi_D:P
\rightarrow P/ D $ be the projection onto the leaf space
of the distribution $D$ and
$\overline{T}_2 $ the  symbol  defined in \eqref{t2 over
line}. Suppose  that $(P, D)$  is
$\overline{T}_2$-idempotent. Then:
\begin{description}
\item[(i)] The relation $ \mathcal{R}_{\overline{T}_2} $ on $P$ 
defined by
$x \mathcal{R}_{\overline{T}_2} y$ if and only if $y \in
\overline{T}_2(x)$ is an equivalence relation.
\item[(ii)]  The following properties are equivalent:
\begin{enumerate}
\item  $y \notin \overline{T}_2(x)$.
\item  $\overline{T}_2(x) \not= \overline{T}_2(y)$.
\item  $\overline{T}_2(x) \cap \overline{T}_2(y)=\emptyset$.
\item  There exist  open neighborhoods $V_x,V_y$ of $x$ and
$y$, respectively, such that $\overline{T}_2(V_x) \cap
\overline{T}_2(V_y)=\emptyset$.
\end{enumerate}
\item[(iii)]  If the projection $\pi_{\overline{T}_2}:P
\longrightarrow P/\mathcal{R}_{\overline{T}_2}$ is an open map then
the quotient space $P/\mathcal{R}_{\overline{T}_2}$ is a
Hausdorff topological space.
\end{description}

\end{lemma}
\noindent
\textbf{Proof}
{\bf (i)} Only transitivity needs to be proved. Let $x,y,z \in P $
be such that  
$x
\mathcal{R}_{\overline{T}_2} y$ and $y \mathcal{R}_{\overline{T}_2} 
z$. By definition,
$ \pi_{D}(x) \mathcal{R}_{T_2} \pi_{D}(y)$ and
$ \pi_{D}(y) \mathcal{R}_{T_2} \pi_{D}(z)$. 
Since the $\overline{T} _2 $-idempotency of $(P, D)$ is equivalent to
the
$T_2$-idempotency of $P/ D  $, Lemma~\ref{equivalence
relation} guarantees that
$
\pi_{D}(x) \mathcal{R}_{T_2} \pi_{D}(z)$ and hence $
\pi_{D}(x) \in T_2( \pi_{D}(z))$. Consequently,  $x
\in 
\pi^{-1}_D\left( T_2 (\pi_{D}(z))\right)=\overline{T} _2 (z)$ and
thus
$z
\mathcal{R}_{\overline{T}_2} x$. 

In order to prove parts {\bf (ii)} and {\bf (iii) }  it suffices
to mimic the corresponding implications in Lemma~\ref{equivalence
relation} but, this time, keeping in mind that the projection
$\pi_{\overline{T}_2}:P
\rightarrow P/\mathcal{R}_{\overline{T}_2}$,
$\pi_{\overline{T}_2}=\pi_{T_2}\circ
\pi_{D}$, is just the composition of two projection maps
and that $\pi_{D}$ is an open map. $\blacksquare$

\begin{theorem}
\label{casimires t2 therorem} 
Let $P$ be a smooth Hausdorff paracompact finite
dimensional manifold and $D$  a smooth integrable
generalized distribution on $P$. Let $\pi_D:P
\rightarrow P/ D $ be the projection onto the leaf space
of the distribution $D$ and
$\overline{T}_2 $ the  symbol  defined in \eqref{t2 over
line}. Suppose  that $(P, D)$  is
$\overline{T}_2$-idempotent and that $\pi_{T_2}$ (and hence
$\pi_{\overline{T}_2}$) is open. Then the continuous 
first integrals of
$D$ separate the
$\overline{T}_2$ sets. In this situation, for any $z \in P $, there
exist continuous first integrals $\{ C _i\}_{i \in  I }
\subset C ^0(P) $ of $D$ such that 
\begin{equation}
\label{we localized end}
\overline{T}_2 (z)=\bigcap _{i \in I}C _i^{-1}(C _i (z)).
\end{equation}
\end{theorem}
\noindent
\textbf{Proof} Since $P$ is by hypothesis
paracompact, so are the quotient spaces $P/D$ and
$P/\mathcal{R}_{\overline{T}_2}$. The hypothesis on the
$\overline{T}_2$-idempotency of $(P, D)$ implies, by
Lemma~\ref{lemma upstairs}, that the  quotient space $
P/\mathcal{R}_{\overline{T}_2} $  is also Hausdorff. Since a
Hausdorff paracompact space is normal, Urysohn's Lemma guarantees
the existence of continuous functions
$f$ on $P/\mathcal{R}_{\overline{T}_2}$ that  separate
its points. The pull back $f \circ
\pi_{\overline{T}_2}\in C ^0(P)$ is a first integral of
$D$. The family of functions of the form $ f \circ
\pi_{\overline{T}_2} $ where $f :
P/\mathcal{R}_{\overline{T}_2} \rightarrow \mathbb{R}$
is a continuous function that separates
two arbitrary points, is the family of continuous first
integrals of $D$ in the statement of the theorem.

In order to prove the identity~(\ref{we localized end}) it
suffices to reproduce the proof of Lemma~\ref{lemma famous inclusion
for ever}, taking   this time   the function $C:P
\rightarrow 
\mathbb{R}^I
$ whose components are the continuous first integrals of $D$ that
separate the $\overline{T}_2$ sets and whose existence we just
proved. 
$\blacksquare$

\begin{remark}
\normalfont
The two hypotheses in the statement of this result, that is the
$\overline{T}_2  $-idempotency and the openness of the projection $
\pi _{\overline{T}_2 }$  are
independent. Indeed, consider the foliation of the Euclidean plane
$\mathbb{R} ^2  $ by the integral curves of the vector field
$\phi (x) \partial/\partial x $, where $\phi$ is a smooth function 
satisfying  $\phi(x)=0$, for  $x\leq 0$, and
$\phi(x)>0$, for  $x>0$. 
In this situation
$\overline{T}_2(x,y )=\{(x,y)\}$, when $x<0$, and
$\overline{T}_2(x,y)=\{(x,y) \in  \mathbb{R} ^2\mid x\leq 0\}$,
if
$x
\geq 0$. In this situation, we obviously have
$\overline{T}_2$-idempotency. However,  the projection $
\pi _{\overline{T}_2 }:\mathbb{R}^2
\rightarrow \mathbb{R} ^2/
\mathcal{R}_{\overline{T}_2}
$  is
not open. Indeed, the saturation $ \pi
_{\overline{T}_2 } ^{-1}(\pi
_{\overline{T}_2 } (U))=\{(x,y)\in  \mathbb{R} ^2\mid
x\geq 0\}$ of the open set $U=\{(x,y )\in  \mathbb{R}^2
\mid x>0\}$ is closed, which is not compatible with $\pi
_{\overline{T}_2} $ being open.
\end{remark}

The following result provides another sufficient condition for the
conclusion of Theorem~\ref{casimires t2 therorem} to hold.

\begin{theorem}
\label{extension theorem} 
Let $D$ be a generalized smooth integrable distribution defined on
the second countable finite dimensional manifold $P$. 
Suppose that there
exist
continuous first integrals $C_i \in C^{0}(P)$, $ i \in  I $, of the
foliation induced by
$D$ that separate its regular leaves. Additionally, assume 
that the map
$C:P
\rightarrow \mathbb{R} ^I$ defined by $C (z):=(C_i(z))_{i \in  I} $,
$z 
\in P$, is open onto its image when $\mathbb{R} ^I  $ is endowed with
the product topology. Then for any $z \in P $
\[
\overline{T}_2 (z)=\bigcap _{i \in I}C _i^{-1}(C _i (z)).
\]
\end{theorem}
 
\noindent
\textbf{Proof} Notice first that the inclusion 
\[
\overline{T}_2 (z)\subset\bigcap _{i \in I}C _i^{-1}(C _i (z)).
\]
is a particular case of~(\ref{famous inclusion for ever}). 

In order to prove the converse inclusion let 
$\pi_{D}: P \rightarrow P/ D$ be
the projection onto the leaf space and
$c:P/ D \rightarrow \mathbb{R} ^k $ the
continuous mapping uniquely determined by the relation
$c \circ \pi_{D}=C $. Let $n \in 
\bigcap_{i\in  I}  C^{-1}_i(C _i(z)) $ , that is, $C (n)= C
(z) $ and assume that $n \notin
\overline{T}_2(z) $. This implies the existence of
two open neighborhoods $V_{\pi_{D} (n)} $  and
$V_{\pi_{D} (z)} $ of  $\pi_{D} (n) $
and $\pi_{D} (z)$, respectively, such that 
$V_{\pi_{D} (n)} \cap V_{\pi_{D} (z)}
=\emptyset $. We will assume for the time being that the
leaf $\pi_{D} (n) $  is regular and will prove
that the assumption $n \notin
\overline{T}_2(z) $ leads to a contradiction. We
will prove later on that the  situation in which 
$\pi_{D} (n) $ is a singular
leaf can be reduced to this case.

If $\pi_{D} (n) $  is regular, the set 
$V_{\pi_{D} (n)}^{reg} $ of regular leaves in
$V_{\pi_{D} (n)}$ is an open dense neighborhood
of $\pi_{D} (n) $ in $V_{\pi_{D} (n)}$.
The openness hypothesis on the map $C$ implies that the
set
\[
U _{C (z)}:= c (V_{\pi_{D}
(n)}^{reg} )\cap  c (V_{\pi_{D} (z)})
\]
is an open neighborhood of $C (z) $. Moreover, the
continuity of
$c$  implies that the sets
\[
A:= c ^{-1}(U _{C (z)})\cap V_{\pi_{D}
(n)}^{reg}\quad\text{and} \quad B:= c ^{-1}(U _{C(z)})
\cap V_{\pi_{D} (z)}
\]
are open neighborhoods of $\pi_{D}
(n)$ and $\pi_{D}
(z)$, respectively. Let $\pi_{D}
(z')$ be  a regular leaf in $B$. The construction of $B$
implies that there exists a regular leaf $\pi_{D}
(s) \in V_{\pi_{D}
(n)}^{reg} \subset V_{\pi_{D}
(n)}  $ such that  $c(\pi_{D}
(z'))=c(\pi_{D}
(s))$. The separation hypothesis on the map $C$ implies 
that $\pi_{D}
(s)=\pi_{D}
(z') \in V_{\pi_{D} (n)} \cap V_{\pi_{D} (z)} $ which is a  contradiction.

In order to conclude the proof we need to show that the
case in which $\pi_{D} (n)$ is singular can be reduced to the situation
that we just treated. Indeed, take $U _{C (z)}:= c
(V_{\pi_{D} (n)} )\cap  c (V_{\pi_{D}
(z)}) $. By the openness of $C$,
$U _{C (z)} $ is an open neighborhood of $C (z) $.
Additionally,  the continuity of
$c$  implies that the sets
$A:= c ^{-1}(U _{C (z)})\cap V_{\pi_{D}
(n)}$ and $ B:= c ^{-1}(U _{C (z)})\cap
V_{\pi_{D} (z)}$ are open disjoint neighborhoods
of $\pi_{D} (n)$ and $\pi_{D}
(z)$, respectively. Let $\pi_{D}
(z')$ be  a regular leaf in $A$. The construction of $A$
implies the existence of a leaf $\pi_{D}
(s) \in V_{\pi_{D} (z)}  $  such that  $C (\pi_{D}
(z'))= C(\pi_{D}
(s))$. If we follow the  preceding argument, replacing
$\pi_{D} (z') $ by $\pi_{D}
(n) $, $\pi_{D}
(s)$ by $\pi_{D}
(z)$, $A$ by $V_{\pi_{D} (n)}$, and with
$V_{\pi_{D} (z)}$ playing the same role we also
obtain a contradiction with the hypothesis
$V_{\pi_{D} (n)} \cap V_{\pi_{D} (z)}
=\emptyset $. $\blacksquare $

\medskip

The reader may be wondering how the two sufficient conditions
for~(\ref{we localized end}) to hold that we presented in the
statements of the theorems~\ref{casimires t2 therorem}
and~\ref{extension theorem}  are related. Our next result answers
this question.

\begin{proposition}
\label{what kind of equivalence}
Let $P$ be a smooth second countable  finite
dimensional manifold and $D$  a smooth integrable
generalized distribution on $P$. Suppose that there
exist
continuous first integrals $C_i \in C^{0}(P)$, $ i \in  I $, of the
foliation induced by
$D$ that separate its regular leaves such
that the map
$C:P
\rightarrow C (P) \subset \mathbb{R} ^I$ defined by $C
(z):=(C_i(z))_{i
\in  I}
$,
$z 
\in P$, is open onto its image when $\mathbb{R} ^I  $ is endowed with
the product topology. Then $(P, D)$  is
$\overline{T}_2$-idempotent and
$\pi_{\overline{T}_2}:P \rightarrow  P/\mathcal{R}_{\overline{T}_2}$
is an open map.
\end{proposition}

\noindent\textbf{Proof.\ \ } In the hypotheses of the statement,
Theorem~\ref{extension theorem} implies that 
$
\overline{T}_2 (z)=C ^{-1}(C (z))$, for any $z \in P $. In particular
\[
\overline{T}_2 (\overline{T}_2(z))= \overline{T} _2(C ^{-1}(C
(z)))=\bigcup_{y \in C ^{-1}(C (z))}\overline{T}_2 (y)=C ^{-1}(C
(z))=\overline{T}_2(z),
\] 
which guarantees that $(P, D) $ is $\overline{T}_2$-idempotent and
hence allows us to define an equivalence relation
$\mathcal{R}_{\overline{T}_2} $ on $P$. We will now show that the
associated projection to the quotient $\pi_{\overline{T}_2}:P
\rightarrow  P/\mathcal{R}_{\overline{T}_2}$ is open. Let $\varphi:
P/\mathcal{R}_{\overline{T}_2} \rightarrow  C (P) $ be the map
defined by  $\varphi(\pi_{\overline{T}_2}(z)):=C (z)  $, $z \in P $.
The equality
$
\overline{T}_2 (z)=C ^{-1}(C (z))$, $ z \in  P $, guarantees that
$\varphi  $ is  a well defined bijection that makes the diagram
\unitlength=5mm
\begin{center}
\begin{picture}(9,6)
\put(1,5){\makebox(0,0){$P$}} 
\put(9,5){\makebox(0,0){$C(P)$}}
\put(5 ,0){\makebox(0,0){$P/\mathcal{R}_{\overline{T}_2}$}} 
\put(1,4){\vector(1,-1){3}}
\put(6,1){\vector(1,1){3}} 
\put(2,5){\vector(1,0){6}}
\put(1,2){\makebox(0,0){$\pi_{\overline{T}_2}$}}
\put(9,2){\makebox(0,0){$\varphi$}}
\put(5,6){\makebox(0,0){$C $}}
\end{picture}
\end{center}
commutative. The continuity and the openness of $C$ imply
respectively the continuity and the openness of $\varphi$, that is,
$\varphi $ is a homeomorphism. Since $\pi_{\overline{T}_2}= \varphi
^{-1} \circ C $, the openness of  $\pi_{\overline{T}_2} $ follows.
\quad $\blacksquare$

\begin{remark}
\normalfont
The converse of the implication in the previous proposition is not
true in general. A counterexample to this effect is an irrational
foliation of the two-torus. In that particular case the $\overline{T}
_2 $ set of any point is the entire torus and hence we have
$\overline{T} _2 $--idempotency with a projection
$\pi_{\overline{T}_2}:P \rightarrow  P/\mathcal{R}_{\overline{T}_2}$
that is obviously open. Nevertheless, the only first integrals of
this foliation are the constant functions that do not separate the
leaves of the foliation, all of which happen to be regular in this
case.
\end{remark}

We now collect the results in Theorems~\ref{casimires t2 therorem}
and~\ref{extension theorem}  and in Proposition~\ref{what kind of
equivalence} and we apply them to the situation in which $P$  is a
Poisson manifold foliated by its symplectic leaves. 
The following result provides two sufficient conditions for the
continuous and topological energy-Casimir methods to coincide.

\begin{theorem}
\label{equivalence result for Poisson}
Let $(P,\{ \cdot , \cdot \})$ be a Poisson (paracompact, second
countable, and Hausdorff) manifold. Let   $\overline{T}_2$ be 
the symbol associated to the symplectic foliation of $P$
induced by the  Poisson structure
$ \{
\cdot , \cdot \} $. 
\begin{description}
\item [(i)] Suppose that there
exist
continuous Casimir functions $C_i \in C^{0}(P)$, $ i \in  I $, that
separate the  regular symplectic leaves of  $P$ such that the map
$C:P
\rightarrow C (P) \subset \mathbb{R} ^I$ defined by $C
(z):=(C_i(z))_{i
\in  I}
$,
$z 
\in P$, is open onto its image when $\mathbb{R} ^I  $ is endowed with
the product topology. Then $P$  is
$\overline{T}_2$-idempotent and
$\pi_{\overline{T}_2}:P \rightarrow  P/\mathcal{R}_{\overline{T}_2}$
is an open map.
\item [(ii)] If $(P,\{ \cdot , \cdot \})$  is
$\overline{T}_2$-idempotent and
$\pi_{\overline{T}_2}:P \rightarrow  P/\mathcal{R}_{\overline{T}_2}$
is an open map then there exist continuous Casimir functions $\{ C
_i\}_{i \in  I }
\subset C ^0(P) $ of $(P, \{ \cdot , \cdot \})$ such that for any $z
\in P $
\[
\overline{T}_2 (z)=\bigcap _{i \in I}C _i^{-1}(C _i
(z)).
\]
\end{description}
\end{theorem}

\begin{remark}
\normalfont
As one could expect, the hypotheses of this theorem are not satisfied
by Montaldi's example (see Figure~\ref{fig:james
example}). Indeed, in this particular case  $\overline{T}
^U_2(\overline{T}
^U_2(0,0,0))= U\neq \overline{T} ^U_2(0,0,0)
$, for any open neighborhood $U$ of the origin $(0,0,0)$. 
\end{remark}

\begin{example}
\normalfont 
The following example is given in~\cite{topological stability}
and provides a situation where their topological 
energy-Casimir method (Theorem~\ref{patrick theorem statement})
works when establishing the stability of a Poisson
equilibrium while the standard energy-Casimir tool fails. We will
show that in this situation the generalized energy-Casimir method
formulated in Corollary~\ref{how to compare both} also works
and, moreover, both results can be applied
interchangeably since  the hypotheses of part {\bf (i)} in 
Theorem~\ref{equivalence result for Poisson} are
satisfied.  

Let
$(\mathbb{R}^3,
\{
\cdot ,
\cdot \}, h)$ be the Poisson dynamical system given by 
\[
\{f,g\}=\nabla A \cdot(\nabla f \times \nabla g),\qquad
A(x,y,z)=(a^2 x^2 -y^2) y
\]where $a$ is a nonzero real constant and $h(x,y,z)=x^2 - y^2 +
z^2$.
Notice that the function $A$
is a Casimir of the bracket $\{ \cdot , \cdot \} $ and that the
points on the form $(x,0,z)$  are equilibria of the Hamiltonian
vector field $X _h  $. We will focus on the stability of the origin
$(0,0,0) $ that happens to be a singular point of the
symplectic foliation of $\mathbb{R}^3 $. In order to
verify that the hypotheses of Corollary~\ref{how to
compare both} are satisfied notice that the map $A$  can
be rewritten as
$A(x,y,z)=(ax+y)(ax-y)y$ and hence its zero level set (the one
containing the equilibrium $(0,0,0)$) can be written as
the union of three irreducible algebraic varieties 
$V_1$, $V_2$, and $V_3$ that are the zero level sets of the
functions
$y, ax-y$, and $ax+y$, respectively. Consequently,
\begin{equation}
\label{decomposition algebraic}
h ^{-1} (0)\cap A ^{-1} (0)=h ^{-1} (0)\bigcap \left(V _1\cup V
_2\cup V _3
\right)=(h^{-1}(0)
\cap V_1) \cup (h^{-1}(0) \cap V_2) \cup
(h^{-1}(0) \cap V_3) 
\end{equation}
which is a single point whenever $|a|<1$ hence
proving the Lyapunov stability of
$(0,0,0)$. 
This is so since each of the three 
intersections on the right hand side of
expression~(\ref{decomposition algebraic}) coincide with the point
$m$. This statement can be proved by showing that the
Hamiltonian restricted to the submanifolds $V _1, V _2$  and $V _3
$ has  a non degenerate critical point at $(0,0,0)$.
This is closely related to the {\bfi  smoothing} of the
$\overline{T} _2 $ set introduced in~\cite{topological
stability}. 

Since the Casimir function $A$
clearly separates the regular symplectic leaves of $(\mathbb{R}
^3,\{ \cdot , \cdot \} )$ and it is an open map, by
Theorem~\ref{equivalence result for Poisson} we can conclude 
that 
\[
\overline{T}_2(0,0,0)=A^{-1}(A(0,0,0))
\] 
and hence energy-Casimir and $T _2 $-based sufficient stability
conditions can be used interchangeably. 
\end{example}

\bigskip

\addcontentsline{toc}{section}{Acknowledgments}
\noindent\textbf{Acknowledgments.}  We thank George Patrick for his
advice and for carefully reading the first draft of this paper. His
suggestions have greatly improved this work. We thank Jerry
Marsden, James Montaldi,  Mark Roberts, and Claudia Wulff for many
illuminating  discussions about these topics over the years. This
research was partially supported by the European Commission through
funding for the Research Training Network
\emph{Mechanics and Symmetry in Europe} (MASIE) and by the Marie Curie 
fellowship
HPTM-CT-2001-00233.

\end{document}